\documentclass[10pt, leqno, twoside]{article}
\usepackage{amssymb, amsmath} 
\usepackage{latexsym} 
\usepackage{graphics} 
\usepackage{url} 
\usepackage{enumerate}
\usepackage{tikz-cd} 

\usepackage[all]{xy}

\usepackage{graphicx}

\usepackage{float}

\usepackage{amsrefs}

\usepackage{color}

\usepackage{bm}

\usepackage{comment}

\usepackage{hyperref}

\numberwithin{equation}{section}

\newtheorem{theo}{Theorem}[section]
\newtheorem{coro}[theo]{Corollary}

\newtheorem{lemm}[theo]{Lemma}
\newtheorem{prop}[theo]{Proposition}
\newtheorem{exam}[theo]{Example}

\newenvironment{proo}[1][\proofname]{\normalfont{\itshape
#1{:}}\quad\mdseries\ignorespaces}
{{\hspace{13cm}$\Box$}{\vskip\belowdisplayskip}}
\newcommand{\proofname}{Proof}

\newtheorem{defi}[theo]{Definition}

\newtheorem{rem}[theo]{Remark}

\newcommand{\N}{{\Bbb N}}

\newcommand{\be}{{\beta}}

\begin{document}

\title{Riordan patterns's quest within simplicial complexes}

\author{Pedro J. Chocano, Ana Luz\'on, Manuel A. Mor\'on and Luis Felipe Prieto-Mart\'inez}

\date{}
\maketitle

\begin{abstract} 

The aim of this paper is twofold. First, we demonstrate how Riordan matrices can be employed to connect well-known concepts in geometric combinatorics, such as $f$-vectors, $h$-vectors, $\gamma$-vectors, in a similar fashion to the McMullen Correspondence,
and the Dehn–Sommerville equations, among others. Second, we investigate the combinatorial properties of the topological join operation, both for simplicial complexes and for Alexandroff spaces. Finally, we explore the Riordan matrices arising from the iteration of this topological operation and analyze their properties.

\end{abstract}

\section{Introduction}

This paper aims to explore the connection between certain combinatorial properties of finite simplicial complexes and Riordan matrices. It includes a section (Section \ref{sect.intro}) that describes some basic facts or preliminaries, followed by two main parts (Section \ref{sect.re} and Sections \ref{sect.properties} and \ref{sect.qc}) that present original results. In this section, we outline the main results, providing both motivation and context.

Since we need to treat many different type of objects, let us include now some general principles concerning notation, to try to avoid confusing the reader:

\begin{itemize}

\item $\mathcal F$ always denote a simplicial complex.

\item $X,Y$ always denote Alexandroff spaces.

\item $P, P_0,P_1,P_2,\ldots$ always denote polynomials. The numerical subindex, if any, is the degree.

\item $P_{\mathcal F},P_X$ always denote the $f$-polynomial of $\mathcal F,X$, respectively.

\item $\widetilde P_{\mathcal F},\widetilde P_X$ always denote the extended $f$-polynomial of $\mathcal F,X$, respectively.

\item $Q_{\mathcal F}$ always denote the $h$-polynomial of $\mathcal F$.

\item $\Gamma_{\mathcal F}$ always denote the $\gamma$-polynomial of $\mathcal F$.

\item Bold letters are reserved for vectors and matrices.

\item $\Phi,\Psi$ always denote functions.

\item $\alpha,\zeta,\kappa,\tau,\omega$ always denote formal power series.

\item $\mathbb K$ always denote a field of characteristic 0 and $\mathbb N$ always denote the set of natural numbers including 0.

\end{itemize}

\subsection{First part}



There is a broad problem in the field of Combinatorial Topology, which can be stated as follows: \emph{Given a family of simplicial complexes (typically those satisfying some fixed topological or combinatorial property), determine the set of integer sequences corresponding to the $f$-vectors of this family, or at least provide necessary conditions they must satisfy.} This problem has been extensively studied in the literature, for various families of simplicial complexes (see \cite{BE, B3, B2, B1, BK, P, Z}, for instance). For some significant families of simplicial complexes, such as the entire set of finite simplicial complexes or the set of simplicial complexes with a prescribed sequence of Betti numbers, the answer is known (see Kruskal-Katona Theorem in \cite{Z} and Bj\"orner-Kalai Theorem in \cite{BK}, respectively).

When studying this problem, the conditions for the $f$-vector are stated in different ways. In this part, we demonstrate that some of the most well-known conditions for the $f$-vector in this context can be expressed in terms of Riordan matrices. First, we explain how the $h$-vector, $g$-vector, and $\gamma$-vector can be derived from the $f$-vector by left-multiplying it by a (finite or infinite) Riordan matrix (Theorem \ref{theo.changes}). This process is similar to the \emph{McMullen correspondence}, which, in some contexts, assigns an $f$-vector to each $g$-vector, also by multiplying the latter by a matrix. Second, we express the well-known \emph{Dehn-Sommerville Equations} (in terms of $f$-vectors, see Statement (a) in Definition \ref{defi.ds}) as an eigenvector problem for a finite Riordan matrix, which is also an involution (Lemma \ref{lemm.dsr}). We then use this fact to study some aspects of a problem previously examined by M. Grunbaum in \cite[Section 9.5]{Gr}, concerning the linear span of the $f$-vectors satisfying the \emph{Dehn-Sommerville Equations} for a fixed dimension (Theorem \ref{theo.gr}). Third, we explore the possibility of using Theorem \ref{theo.changes} and the Riordan matrix framework to study some problems and results discussed in the literature (Remarks \ref{rem.RMcM} and \ref{rem.pos}).

\subsection{Second part}

We explore the operation that we call the \emph{$q$-cone}. This operation generalizes the well-known operations of \emph{cone} and \emph{suspension}, which are defined for finite simplicial complexes. The $q$-cone operation is introduced in Section \ref{sect.properties} (Definition \ref{defi.qcone}). We begin by proving that the$f$-vector of the $q$-cone of a simplicial complex $\mathcal F$ can be easily described in terms of the $f$-vector of $\mathcal F$ (Lemma \ref{L:f-q-cones}). Similarly, we show the analogue for the homology groups of the geometric realizations (Lemma \ref{thm_homologia_conos}). To enhance  understanding of this operation and because it illustrates  some concepts of Combinatorial Topology, we also prove that the $q$-cone preserves partitionability for pure simplicial complexes (Theorem \ref{thm_join_particionable}) and vertex-decomposability (Theorem \ref{thm_vertex_decomposable_join}). This, in turn, produces a result that allows us to describe the $h$-vector of the $q$-cone of  $\mathcal F$ in terms of  the $h$-vector of $\mathcal F$ (Theorem \ref{thm_exact_facing_h_vector_join}).

In Section \ref{sect.qc}, for any two natural numbers $m$ and $q$, we describe a sequence $(\Delta_{m,q}^{(n)})$ of simplicial complexes of increasing dimension, defined by iteratively performing the $q$-cone operation (Definition \ref{defi.qsec}). This family of simplicial complexes serves as a well-understood example that helps us grasp many concepts from the literature on simplicial complexes. We show that if we construct a matrix $\bm F_{m,q}$ whose $i$-row corresponds to the $f$-vectors of $\Delta_{m,q}^{(i)}$ for $i\in\mathbb N$, then the resulting matrix is a Riordan matrix. This fact allows us to describe in detail the $f$-vectors of these simplicial complexes $(\Delta_{m,q}^{(n)})$ (Theorem \ref{theo.muchos}). Additionally, the matrix $\widetilde{\bm F}_{m,q}$, obtained by suitably extending $\bm F_{m,q}$ and whose $(i+1)$-row corresponds to the extended $f$-vectors of $\Delta_{m,q}^{(i)}$, for $i\in\mathbb N$, is also a Riordan matrix (Theorem \ref{theo.muchos}). In a similar fashion, if we construct the matrix $\bm H_{m,q}$ using, this time, the $h$-vectors (placed in the adequate way) we again obtain a Riordan matrix. We prove that, surprisingly, for $m,q\geq 2$, $\bm H_{m,q}=\widetilde{\bm F}_{m-1,q-1}$ (Theorem \ref{theo.hmq}). As we will discuss later, this result is of particular interest because it provides a family of examples of simplicial complex whose $h$-vector is also an $f$-vector.

The matrix $\widetilde{\bm F}_{m,q}$ allows us to easily compute the Euler characteristic of all the elements in the sequence $(\Delta_{m,q}^{(n)})$ simultaneously, using the so called \emph{Fundamental Theorem of Riordan Matrices} (this is done in Theorem \ref{theo.vectorChi}). At this point, we find a curious fact: there is a sequence such that, when left multiplied by $\widetilde{\bm F}_{m,q}$, produces the infinite vector $[1,1,\ldots]^T$, that is, when left multiplied by any (row) $f$-vector of an element in $(\Delta_{m,q}^{(n)})$, the same number is obtained (Corollary \ref{coro.meuler}).

Regarding other topological properties of the elements in the sequence $(\Delta_{m,q}^{(n)})$, if we define a fourth matrix $\bm{B}_{m,q}$, whose $i$-th row corresponds to the sequence of reduced Betti numbers of $(\Delta_{m,q}^{(i)})$, we again obtain a Riordan matrix, which this time is diagonal (Theorem \ref{theo.matrixB}).

In Subsection \ref{subsect.qalexandroff}, we give another approach to some of the Riordan matrices obtained before within the combinatorial context of Alexandroff spaces. Specifically, by using the combinatorial analogue of the join (the non-Hausdoff join), we construct spaces (also denoted by) $\Delta_{m,q}^{(n)}$ satisfying that their order complexes are precisely the simplicial complexes defined in Section \ref{sect.qc}. We analyze topological properties of these spaces (homotopy type, weak homotopy type, group of homeomorphisms, the group of self-homotopy equivalences and the determinant) and combinatorial properties (order dimension). From this study, we derive other Riordan matrices by using some of the invariants obtained through this section. 

All the results and facts in Section \ref{sect.qc} are illustrated in Subsection \ref{subsect.exam}, for the cases $m=q=1$ (simplices), $m=q=2$ (triangulated spheres) and $m=1,q\geq 1$.


\section{Basic Facts} \label{sect.intro}

\subsection{Simplicial complexes and $f$-vectors}\label{sect.intro.simplicial.complex}

Let us assume that the reader is familiar to the concepts of \emph{finite simplicial complex} (both, \emph{abstract} and \emph{geometric} and their interaction), \emph{face}, \emph{dimension of a face}, \emph{dimension of a finite simplicial complex}, \emph{pure simplicial complex}, \emph{vertex}, \emph{facet}, \emph{geometric realization}, \emph{subdivision} and \emph{isomorphism of simplicial complexes}. In other case, we refer the reader, for example, to any of the texts \cite{H,ST,S,Z}.






Almost every simplicial complex in this paper is finite so, as a general notational principle, let us avoid including the term \emph{finite}, if possible. The only appearences of a non-finite simplicial complexes will be clearly indicated. 
In the following, when writing \emph{simplicial complex}, we refer to either an abstract or geometric finite simplicial complex. We will specify the nature of the object, if neccessary. Usually, we use the letter $\mathcal F$ to denote a finite simplicial complex and $|\mathcal F|$ for its geometric realization. 




Let us introduce now the most important combinatorial invariant for finite simplicial complexes appearing in this paper:

\begin{defi}\label{def.f-vecto.f-pol.etc} Given a simplicial complex $\mathcal F$ (either abstract or geometric) of dimension $d$, for every $k=0,\ldots, d$, let us denote by $f_k$ to the number of faces of dimension $k$. We define:

\begin{itemize} 

\item The \textbf{face-vector}, or simply \textbf{f-vector}, of $\mathcal F$ to be the finite sequence $\bm f=(f_0,f_1,\ldots,f_d)$.

Sometimes we call $f$-vector, instead, to the infinite sequence $\bm f=(f_0,f_1,\ldots,f_d,0,0,\ldots)$. 

\item The \textbf{f-polynomial} of $\mathcal F$ to be the polynomial
$$P_{\mathcal F}(x)=f_0+f_1x+f_2x^2+\ldots+f_dx^d. $$

\item The \textbf{extended f-vector} to be the sequence $(f_{-1},f_0,f_1,\ldots, f_d)$, where $f_{-1}=1$  (representing the number of empty faces in $\mathcal F$ because empty faces are usually considered to have dimension $-1$).

\item The \textbf{extended f-polynomial} of $\mathcal F$ to be the polynomial
$$\widetilde P_{\mathcal F}(x)=f_{-1}+f_0x+f_1x^2+\ldots+f_dx^{d+1}=1+xP_{\mathcal F}(x). $$

\end{itemize}

\end{defi}

Note that the degree of the f-polynomial coincides with the dimension of $\mathcal F$ and that the degree of the extended f-polynomial is the dimension of $\mathcal F$ plus 1.



\subsection{Alexandroff spaces and posets}\label{subsection.preliminares.posets}

We recall the basic theory and definitions of Alexandroff spaces (for a general introduction we recommend \cite{Barmak,May}). Let $X$ be a topological space. It is said that a topological space $X$ is an \textbf{Alexandroff space} whenever the arbitrary intersection of open sets is again an open set. Note that the most important example of Alexandroff spaces are finite topological spaces. From now on, we assume that every Alexandroff space throughout this manuscript satisfies the $T_0$-separation axiom. This condition is crucial in what follows. 

For any Alexandroff space $X$ and $x \in X$, denote $U_x$ as the intersection of every open set containing $x$, which is an open set by hypothesis. Similarly, $L_x$ is defined as the intersection of every closed set containing $x$. Declare $x \leq y$ if and only if $U_x \subseteq U_y$. It is immediate to verify that $X$ with this relation is a partially ordered set. Moreover, for any partially ordered set $(X, \leq)$, the set of upper sets forms a basis for a $T_0$-topology that makes $X$ an Alexandroff space. These two relations are mutually inverse. Moreover, a map $f:X\rightarrow Y$ is a continuous map between Alexandroff spaces if and only if $f$ is an order-preserving map. From this, it can be deduced that the category of Alexandroff spaces is isomorphic to the category of partially ordered sets. This allows us to treat the same object from two different points of view. For this reason, we will not distinguish between an Alexandroff space and a partially ordered set. 
Now, we introduce some structures related to Alexandroff spaces.

Let us restrict our attention to finite topological spaces (finite spaces for short). Given a finite space $X$, the \textbf{Hasse diagram} of $X$ is a directed graph, where the vertices are the points of $X$ and there is a directed edge from $x$ to $y$ whenever $x<y$ and there is no $z$ such that $x<z<y$. In subsequent Hasse diagrams we assume an upward orientation and omit the direction of the edges. Note that this definition can be generalized to countable Alexandroff spaces. The \textbf{height} of a finite space $X$, denoted by $\text{ht}(X)$, is defined as one less the length of a maximal chain in $X$. Similarly, the height of a point $x\in X $ is defined as $\textnormal{ht}(x)=\textnormal{ht}(U_x)$. The \textbf{width} of $X$, denoted by $\text{width}(X)$, is the size of a maximal antichain in $X$. 

Although finite spaces or, generally, Alexandroff spaces may seem simple in terms of algebraic topology, the following results show the opposite. For an Alexandroff space $X$, the \textbf{order complex} of $X$, denoted by $\mathcal{K}(X)$, is the simplicial complex whose simplices are the non-empty chains of $X$. Similarly, the \textbf{face poset} of a simplicial complex $\mathcal{F}$, denoted by $\mathcal{X}(\mathcal{F})$, is the Alexandroff space whose points are the simplices of $\mathcal F$ ordered by the subset relation. For any point $u\in|\mathcal{K}(X)|$, there exists a unique open simplex $(v_0,...,v_n)$ given by $v_0<\ldots<v_n$ containing it. Define $\Phi_X:|\mathcal{K}(X)|\rightarrow X$ by $\Phi_X(u)=v_0$. In \cite[Theorem 2]{Mccord}, it is proved that $\Phi_X$ is a weak homotopy equivalence. Similarly, there is also a weak homotopy equivalence $\Phi_{\mathcal{X}(\mathcal{F})}:|\mathcal{F}|\rightarrow \mathcal{X}(\mathcal{F})$ (see \cite[Theorem 2]{Mccord}). These two results show that homotopy groups and homology groups of Alexandroff spaces and simplicial complexes are the same. It is also worth pointing out that $\mathcal{K}(\mathcal{X}(\mathcal{F}))$ corresponds to the usual barycentric subdivision of $\mathcal{F}$ and that the finite barycentric subdivision of an Alexandroff space $X$ is defined by $\mathcal{X}(\mathcal{K}(X))$.

As it was done in Subsection \ref{sect.intro.simplicial.complex}, we can also adapt the notions introduced in Definition \ref{def.f-vecto.f-pol.etc} within this context.

\begin{defi}\label{def.f-vecto.f-pol.etc.finito} Given a finite space $X$ of height $d$, for every $k=0,\ldots, d$, let us denote by $f_k$ to the number of chains of length $k$. We define:

\begin{itemize} 

\item The \textbf{face-vector}, or simply \textbf{f-vector}, to be the finite sequence $(f_0,f_1,\ldots,f_d)$.

Sometimes we also call $f$-vector to the infinite sequence $(f_0,f_1,\ldots,f_d,0,0,\ldots)$. 

\item The \textbf{f-polynomial} to be the polynomial
$$P_{X}(x)=f_0+f_1x+f_2x^2+\ldots+f_dx^d. $$

\item The \textbf{extended f-vector} to be the sequence $(f_{-1},f_0,f_1,\ldots, f_d)$, where $f_{-1}=1$  (representing the number of empty faces in $\mathcal F$ because empty faces are usually considered to have dimension $-1$)

\item The \textbf{extended f-polynomial} of $X$ to be the polynomial
$$\widetilde P_{X}(x)=f_{-1}+f_0x+f_1x^2+\ldots+f_dx^{d+1}=1+xP_{X}(x). $$

\end{itemize}

\end{defi}

\begin{rem}\label{rem.f.vectores.son.iguales.finitos.mccord} Note that we use the same notation in Definition \ref{def.f-vecto.f-pol.etc} and Definition \ref{def.f-vecto.f-pol.etc.finito} because the f-vector of a finite space $X$ is the same as the f-vector of $\mathcal{K}(X)$.
\end{rem}

Given a finite space $X$, we say that $x\in X$ is an up (down) beat point if $L_x\setminus \{x\}$ has a minimum ($U_x\setminus \{x \}$ has a maximum). Removing a beat point does not change the homotopy type of the space (see \cite{Stong}). We say that a finite space $X$ is a \textbf{minimal finite space} if $X$ does not contain beat points. Let $X$ be a minimal finite space and let $\Phi:X\rightarrow X$ be a continuous map. Then $\Phi$ is a homeomorphism if and only if $\Phi$ is a homotopy equivalence. This gives that the group of homeomorphisms $Aut(X)$ is the same as the group of self-homotopy equivalences $\mathcal{E}(X)$ when $X$ is a minimal finite space.  

To conclude this section, we recall a fundamental result regarding to homotopic maps within this framework. 

\begin{theo}\label{thm.homotopy.alexandroff.spaces}
    Let $\Phi,\Psi:X\rightarrow Y$ be continuous maps between Alexandroff spaces. If there is a sequence of continuous maps  $\Phi_{0},\ldots,\Phi_{n}:X\rightarrow Y$ such that $\Phi_0=\Phi$, $\Phi_n=\Psi$ and for every $i=0,...,n-1$, either  $\Phi_i(x)\leq \Phi_{i+1}(x)$ or $\Phi_i(x)\geq \Phi_{i+1}(x)$ for every $x\in X$, then $\Phi$ is homotopic to $\Psi$.
\end{theo}

\subsection{Riordan matrices and the Riordan group}

This subsection is a summary of notation and results that is developed, in detail, in, for example, the papers \cite{CLMPS, 2ways, LMMPS, poly, teo}. The definition of {\it Riordan matrix} and the related concept  of {\it Riordan group} appeared in the foundational paper \cite{Sha91} due to L.W. Shapiro, S. Getu, W. J. Woan, and L. Woodson. The original definition of a Riordan matrix given in \cite{Sha91} is more restrictive than that used currently in the literature and in this paper. 

\begin{defi} A \textbf{Riordan matrix} is now defined as an infinite matrix $\bm D=[d_{i,j}]_{i,j\in\N}$ whose columns are the coefficients of successive terms of a geometric progression in $\mathbb K[[x]]$ where the initial term is a formal power series of order $0$ and the common ratio is a formal power series of order $1$. 
\end{defi}

As a consequence of this definition,  $\bm D$ needs to be lower triangular, that is, 
\begin{equation} \label{eq.im}\bm D=\begin{bmatrix} d_{00}\\ d_{10} & d_{11} \\ d_{20} & d_{21} & d_{22} \\ \vdots & \vdots & \vdots & \ddots \end{bmatrix} \end{equation}

\noindent and its diagonal $(d_{00},d_{11}, d_{22},\ldots)$ needs to be a geometric progression in $\mathbb{K}$. These matrices are always invertible.

In the following, we shall not write the dependence of the elements in $\mathbb K[[x]]$ in the variable $x$, that is, we may write $\alpha$ instead of $\alpha(x)$, and so on.

We denote a Riordan matrix $\bm D$ as  $T(\alpha\mid
\omega)$, for $\alpha(x)=\sum_{k=0}^{\infty}\alpha_k x^k$ and $\omega(x)=\sum_{k=0}^{\infty}\omega_k x^k$ being formal power series in $\mathbb K[[x]]$ with $\alpha(0)\neq0$ and $\omega(0)\neq0$, if the first term of the geometric progression in $\mathbb K[[x]]$ mentioned above is $\displaystyle{\frac{\alpha(x)}{\omega(x)}}$
and the common ratio is  $\displaystyle{\frac{x}{\omega(x)}}$. Consequently, 
$\displaystyle{d_{i,j}=[x^i]\frac{x^j\alpha(x)}{\omega^{j+1}(x)}}$. In this notation, Pascal's triangle is $T(1\mid 1-x)$.


In \cite{Sha91}, one of the main results about Riordan matrices was stated. Currently (extended for the definition of Riordan matrix presented above) many authors call it the
{\it Fundamental Theorem for Riordan Matrices (FTRM)}. 

\begin{theo}[FTRM] Let $\bm D=T(\alpha\mid \omega)$ be a Riordan matrix and let $\zeta(x)=\sum_{k=0}^{\infty}
\zeta_k x^k$ be a power series in $\mathbb K[[x]]$. Consider the column vector ${\bf c}=[\zeta_0, \zeta_1, \zeta_2, \ldots]^{\mathrm{T}}$. Then, the generating function of the matrix product $\bm D{\bf c}$ is a second column vector whose generating function is $\frac{\alpha(x)}{\omega(x)}\zeta(\frac{x}{\omega(x)})$.  
\end{theo}

We may interpret the FTRM as a group action (we will discuss why the set of Riordan matrices is a group in a moment). In this case, we may write $T(\alpha\mid \omega)(\zeta)=\frac{\alpha}{\omega}\zeta(\frac{x}{\omega})$. A proof of this result using a special ultrametric space $(\mathbb K[[x]], d)$ can be found  in \cite[Proposition 19]{teo}.

The \textbf{Riordan group} (i.e., the set of all Riordan matrices with the usual product of matrices), denoted by $\mathcal R(\mathbb K)$  or shortly $\mathcal R$, is a subgroup of the group of invertible infinite lower triangular matrices with the usual product of matrices as the
operation. As a consequence of the FTRM, the product is given by 
\begin{equation} \label{eq.prod} T(\alpha\mid \omega)T(\kappa\mid \tau)=T\left(\alpha\kappa\left(\frac{x}{\omega}\right)\Big| \
\omega\tau\left(\frac{x}{\omega}\right)\right),\end{equation} 

\noindent and the inverse is given by 
\begin{equation}\label{eq.inverses}(T(\alpha\mid \omega))^{-1}\equiv T^{-1}(\alpha\mid \omega)=T\left(\frac{1}{\alpha(\frac{x}{A})}\ \Big|\ A \right),\end{equation}
\noindent where
$\displaystyle{\left(\frac{x}{A}\right)\circ\left(\frac{x}{\omega}\right)=\left(\frac{x}{\omega}\right)\circ\left(\frac{x}{A}\right)=x}$. See \cite[Proposition 7]{2ways} for more details.

The sequence of the coefficients of the previous formal power series, denoted by $A$, is the so-called 
\textbf{A-sequence} of $T(\alpha\mid \omega)$. Obviously, the $A$-sequence of $T(\alpha\mid \omega)$ depends only on the power series $\omega$. Moreover, if $A=\sum_{k\geq0}a_{k}x^{k}$, then
\begin{equation} \label{eq.aseq} \forall i,j\geq 1,\qquad  d_{ij}=\sum_{k=0}^{i-j}a_kd_{i-1,j-1+k}, \end{equation}

\noindent where $(a_0,a_1,a_2,\ldots)$ is the A-sequence of the matrix $T(\alpha\mid \omega)$.

Suppose that a lower triangular infinite matrix, as the one appearing in Equation \eqref{eq.im}, is a Riordan matrix $T(\alpha\mid \omega)$ and consider the sequence of polynomials corresponding to the generating functions of each row, that is, 
$$P_0(x)=d_{00},\qquad P_1(x)=d_{10}+d_{11}x,\qquad P_2(x)=d_{20} +d_{21} x+d_{22}x^2,\qquad \ldots $$

\noindent According to the pattern exhibited in Equation \eqref{eq.aseq}, we have (see \cite[Equation (1)]{poly}) that 
$$P_n(x)=\left(\frac{x-\omega_1}{\omega_0}\right)P_{n-1}(x)-\frac{\omega_2}{\omega_0}P_{n-2}(x)-\ldots-\frac{\omega_n}{\omega_0}P_0(x)+\frac{\alpha_n}{\omega_0} $$

\noindent for $(\alpha_0,\alpha_1,\alpha_2,\ldots)$, $(\omega_0,\omega_1,\omega_2,\ldots)$ being the sequences of coefficients of $\alpha,\omega$, respectively.

We need also to introduce in this subsection the concept and notation related to \emph{finite Riordan matrix}. A lower triangular matrix $\bm D_n$ of size $(n+1)\times (n+1)$ is a \textbf{finite Riordan matrix} if there exists a Riordan matrix $T(\alpha\mid \omega)$ (it is not unique) such that $\bm D$ is its leading principal submatrix. In this case, we use the notation $\bm D_n=T_n(\alpha\mid \omega)$. There is a version of the FTRM, of the formulas in the Equations \eqref{eq.prod} and \eqref{eq.inverses} (for the product and the inverse) and of the formula in Equation \eqref{eq.aseq} (in this case, the A-sequence is a finite sequence). We omit details, that can be found, for example,  in \cite{CLMPS, LMMPS}.

All the Riordan matrices appearing in this paper, from this moment on, correspond to the case $\mathbb K=\mathbb C$.


\section{Riordan matrices and $f$-vectors} \label{sect.re}

As we said in the introduction, the first part of this paper consists in this section and here we explore the idea of stating different conditions for the $f$-vector, well known in the literature, in terms of  Riordan matrices.




\subsection{Riordan correspondences between $f$-vectors, $h$-vectors, $g$-vectors and $\gamma$-vectors}

In the literature, the reader may find that, at some moments, it is more convenient to present the information contained in the face-vector in a different way (we will discuss an example in a moment). Some of these alternatives are the so called $h$-vector, $g$-vector and $\gamma$-vector. The definitions are included below and have been taken from \cite{Z}.

First of all, we have:

\begin{defi}[$h$-vector] \label{defi.h}
For a given extended face-vector $\bm f=(f_{-1},f_0,\ldots, f_d)$, we define the corresponding \textbf{h-vector} as the sequence $\bm h=(h_0,h_1,\ldots, h_{d+1})$ defined, for each $0\leq k\leq d+1$, by
$$h_k=\sum_{i=0}^k(-1)^{k-i}\binom{d+1-i}{d+1-k}f_{i-1}. $$

\noindent It will be also convenient to introduce the \textbf{h-polynomial} as $Q_\mathcal F(x)=h_0+h_1x+\ldots+h_{d+1}x^{d+1}$.

\end{defi}

When the simplicial complex is \emph{partitionable} (we explain this concept below, see Subsection \ref{defi.partitionable}), then the entries in the $h$-vector actually have combinatorial meaning, that is, the entries \emph{are counting something}, so they are always non-negative (see \cite{KO,Z}). 

\medskip

\noindent \textbf{Why is the $h$-vector useful?} For some important classes of simplicial complexes, the $h$-vector is symmetric (equivalently, the h-polynomial is \emph{palindromic} or \emph{self-reciprocal}). This condition is called in the literature \emph{Dehn-Sommerville Equations} and this is one of the example that we announced in which a condition is easier to explain in terms of the $h$-vector than in terms of the $f$-vector. More about this will be said in the following subsection.

Second, we have:

\begin{defi}[$g$-vector] \label{defi.g} For a given h-vector $(h_0,\ldots, h_{d+1})$, we define the corresponding g-vector as the sequence $(g_0,\ldots, g_{d+2})$ defined by setting $g_0=h_0$ and, for each $0\leq k\leq d+2$,
$$g_k=h_k-h_{k-1}, $$

\noindent following the convention $h_{d+2}=0$. We also define the \textbf{$\bm g$-polynomial} as $G_{\mathcal F}(x)=g_0+g_1x+\ldots+g_{d+2}x^{d+2}$.

\end{defi}

\noindent \textbf{Why is the $g$-vector useful?} The g-vector was originally introduced for the problem of the characterization of the face-vector of simplicial polytopes (see \cite{Z} for the definiton). This time, the corresponding condition is more complicated than just \emph{being symmetric}, so we will not include it here (we refer the reader, again, to \cite{Z}). Let us remark that the $h$-vector of such simplicial complexes do satisfy the Dehn-Sommerville Equations and, thus, is symmetric. For this reason, some authors define the $g$-vector as the sequence $(g_0,\ldots,g_{\lfloor (d+1)/2\rfloor})$ ($\lfloor (d+1)/2\rfloor$ denotes the floor of $(d+1)/2$) since, for simplicial polytopes, the entries $(g_{1+\lfloor (d+1)/2\rfloor},\ldots, g_{d+2})$ do not contain any extra information.

\medskip

Finally, we have:

\begin{defi}[$\gamma$-vector] \label{defi.gamma} Let $Q_{\mathcal F}(x)$ be some h-polynomial of degree $d+1$.  We define the $\bm \gamma$\textbf{-series} to be the unique formal power series $\Gamma_{\mathcal F}(x)$ satisfying
$$Q_{\mathcal F}(x)=(1+x)^{d+1}\Gamma_{\mathcal F}\left(\frac{x}{(1+x)^2}\right). $$

\noindent The sequence $(\gamma_0,\gamma_1,\ldots)$ corresponding to the coefficients of $\Gamma_{\mathcal F}(x)$ is called the \textbf{$\bm \gamma$-sequence}. If the entries in the $\bm \gamma$-sequence equal to 0 from the position $\lfloor (d+1)/2\rfloor+1$ on, then $\bm \gamma=(\gamma_0,\ldots,\gamma_{\lfloor (d+1)/2\rfloor})$ is called \textbf{$\bm\gamma$-vector}.

\end{defi}

Note that, if $Q_{\mathcal F}(x)$ is a palindromic polynomial of degree $d+1$, then the $\gamma$-series is a polynomial of degree at most $\lfloor (d+1)/2\rfloor$, that is, we have a $\gamma$-vector. Moreover, if $\bm{h}$ has integral coefficients, then so does $\bm \gamma$.

\medskip

\noindent \textbf{Why is the $\gamma$-vector useful?}  As the $g$-vector is related to the problem of describing the face-vectors of simplicial polytopes, the $\gamma$-vector is related to the problem of describing the face-vectors of  \emph{flag complexes}.

\medskip

In the literature, we may find several formulas relating the entries in the $f$-vector, the $h$-vector, the $g$-vector and the $\gamma$-sequence between them (apart from those that have been included in the previous definitions, see \cite{G,P,Z}). In the following, we propose some unified formulas, changing from one vector to another using Riordan matrices:

\begin{theo}[Riordan correspondences] \label{theo.changes}  Let $\mathcal F$ be a simplicial complex of dimension $d$.

\begin{itemize}

\item[(1)] Let  $\widetilde P_{\mathcal F}(x)$ denote the extended $f$-polynomial  of $\mathcal F$ and let $Q_{\mathcal F}(x)$ denote the h-polynomial of $\mathcal F$. Then
\begin{equation}\label{eq.fh2} Q_{\mathcal F}(x)=T((1-x)^{d+2}\mid 1-x)\widetilde P_{\mathcal F}(x).\end{equation}

\noindent Note that $T((1-x)^{d+2}\mid 1-x)$ depends on the dimension of the simplicial complex.

\item[(2)] Additionally, let $G_{\mathcal F}(x)$ be the g-polynomial of $\mathcal F$. Then
$$G_{\mathcal F}(x)=T(1-x\mid 1)Q_{\mathcal F}(x). $$

\item[(3)] Finally, if $\Gamma_{\mathcal F}(x)$ denotes the $\gamma$-series of $\mathcal F$, then
$$Q_{\mathcal F}(x)=T((1+x)^{d+3}\mid (1+x)^2)\Gamma_{\mathcal F}(x). $$

\end{itemize}

\end{theo}

\begin{proo} For the first two cases, we only need to write explicitly the corresponding matrix and compare the expression with the definition of $h$-vector and $g$-vector (Definitions \ref{defi.h} and \ref{defi.g}):

$$\begin{bmatrix}h_0\\ h_1 \\ \vdots \\ h_{d+1}\\ 0 \\ \vdots \end{bmatrix}=\underbrace{\begin{bmatrix}\binom{d+1}{d+1}\\ -\binom{d+1}{d} & \binom{d}{d} \\ \vdots & \vdots & \ddots \\ (-1)^{d+1}\binom{d+1}{0} & (-1)^{d}\binom{d}{0} & \hdots & \binom{0}{0}\\ \vdots & \vdots & & \vdots & \ddots\end{bmatrix}}_{T((1-x)^{d+2}\mid 1-x)} \begin{bmatrix}f_{-1}\\ f_0 \\ \vdots \\ f_{d}\\ 0\\ \vdots  \end{bmatrix},$$



$$\begin{bmatrix} g_0\\ g_1\\ \vdots \\ g_{d+1} \\ g_{d+2}\\ 0 \\ \vdots \end{bmatrix}=\underbrace{\begin{bmatrix}1 \\ -1 & 1 \\ 0 & -1 & 1 \\ \vdots & \ddots  & \ddots & \ddots \\ 0 &\ldots & 0 & -1 & 1 \\ \vdots & & \vdots & \ddots & \ddots & \ddots \end{bmatrix}}_{T(1-x\mid 1)}\begin{bmatrix} h_0 \\ h_1 \\ \vdots \\ h_{d+1} \\ h_{d+2}\\ 0 \\ \vdots \end{bmatrix}.$$

Finally, the equivalence between Definition \ref{defi.gamma} and the condition appearing in the third statement is an immediate consequence of the FTRM.

\end{proo}

\subsection{Formulation of Dehn-Sommerville Equations in terms of Riordan matrices}

Several important classes of simplicial complexes satisfy the so called \emph{Dehn-Sommerville Equations} (see \cite[Section 8.1]{Z}, for more information).

\begin{defi} \label{defi.ds} Let $\mathcal F$ be a simplicial complex and let $(f_{-1},\ldots, f_d)$, $(h_0,\ldots, h_{d+1})$ be the extended $f$-vector and the $h$-vector, respectively. Then we say that $\mathcal F$ satisfies the \textbf{Dehn-Sommerville Equations} if any of the two following equivalent conditions hold:

\begin{itemize}

\item[(a)] For all $k=0,\ldots, d+1$, $\displaystyle{f_{k-1}=\sum_{i=k}^{d+1}(-1)^{d+1-i}\binom{i}{k}f_{i-1}}$.

\item[(b)] For all $k=0,\ldots, d+1$, $h_k=h_{d+1-k}$.

\end{itemize}

\end{defi}

As we have announced in the previous subsection, this is one of the examples of a condition where the alternative presentation of the information made in the $h$-vector makes things simpler, although it can be more obscure (in the sense that the meaning of the entries of the $h$-vector is less clear than the one of the entries of the $f$-vector).

Riordan matrices provide a good language to express the version of Statement (a) of the Dehn-Sommerville Equations:

\begin{lemm} \label{lemm.dsr} Let $\mathcal F$ be a simplicial complex and let $\bm f=(f_{-1},\ldots, f_d)$ denote its extended $f$-vector. Then $\bm f$ satisfies Dehn-Sommerville Equations if and only if $[f_d,\ldots,f_{-1},0,\ldots]^T$ is an eigenvector of eigenvalue 1 of  
$$T(-(1+x)^{d+2}\mid -(1+x)).$$

\end{lemm}

\begin{proo} To prove the result we only need to write explicitely the corresponding matrix,
$$\begin{bmatrix}f_d\\ f_{d-1}\\ \vdots \\ f_{-1} \end{bmatrix} =\begin{bmatrix}\binom{d+1}{d+1}\\ \binom{d+1}{d} & -\binom{d}{d}\\ \vdots & \vdots & \ddots \\ \binom{d+1}{0} & -\binom{d}{0} & \ldots & (-1)^{d+1}\binom{0}{0} \end{bmatrix}\begin{bmatrix}f_d\\ f_{d-1}\\ \vdots \\ f_{-1} \end{bmatrix},$$

\noindent and compare the expression obtained with the one in Definition \ref{defi.ds}.

\end{proo}


Finding the set of vectors $(f_{-1},\ldots, ,f_d)$  satisfying Dehn-Sommerville Equations (even if they do not correspond to the $f$-vector of any simplicial complex) was a problem already studied by M. Gr\"unbaum in \cite[Section 9.5]{Gr}. Riordan matrices allow us to describe this set:

\begin{theo} \label{theo.gr} Let  $(f_{-1},\ldots, ,f_d)$ be the $f$-vector of a given simplicial complex $\mathcal{F}$ and let $\widetilde P_\mathcal{F}(x)$ be the corresponding extended $f$-polynomial. Then:

\noindent (a) $[f_d,\ldots,f_0,f_{-1}]^T$ is a linear combination of the even columns $0,2,\ldots,\lfloor\frac{d+1}{2} \rfloor$ of the matrix 
$$T_{d+1}\left((1+x)^{\frac{d}{2}+1}\mid\sqrt{1+x}\right).$$

\noindent (b) The polynomial  $x^{d+1}f(1/x)$, equals the Taylor polynomial of degree $d+1$ of the series
$$(1+x)^{(d+1)/2}P(x^2/(1+x)).$$

\noindent  where $P$ is a polinomial of degree $\lceil\frac{d}{2}\rceil$.

\end{theo}

\begin{proo} Note that
$$T_{d+1}(-(1+x)^{d+2}\mid -(1+x))T_{d+1}\left((1+x)^{\frac{d}{2}+1}\mid\sqrt{1+x}\right)=T_{d+1}\left((1+x)^{\frac{d}{2}+1}\mid\sqrt{1+x}\right)T_{d+1}(-1,-1).$$



\end{proo}

\subsection{Possible future applications for the Riordan correspondences: McMullen correspondence and positivity of vectors in terms of Riordan matrices}


Let us begin this subsection by introducing the \emph{McMullen correspondence}. Most of the discussion appearing in this subsubsection can be found in \cite[Section 8.6]{Z}.

One of the most important results in Geometric Combinatorics is, probably, the so called $g$-Theorem. The $g$-Theorem characterizes the $g$-vector of the boundary of a simplicial polytope of dimension $d$ (thus, a simplicial complex of dimension $d-1$): if $(g_0,g_1,\ldots, g_{d+1})$ is the $g$-vector of the boundary of a simplicial polytope of dimension $d-1$, then its \emph{first half}, that is, the sequence $(g_0,\ldots,g_{\lfloor\frac{d}{2}+1\rfloor})$, should be what it is called in the literature an \emph{$M$-sequence} (see  \cite[Theorem 8.35]{Z}).  So we have the:

\medskip

\noindent \textbf{McMullen correspondence.} \emph{There is a bijection }
\begin{equation}\label{eq.McM} \bm g\mapsto \bm g\cdot \bm M_d \end{equation}

\noindent  \emph{(were $\bm M_d$ is a well studied matrix of size $(\lfloor\frac{d}{2}+1\rfloor)\times (d+1)$, see \cite{BE, G2} and the references therein) between the $M$-sequences $\bm g=(g_0,\ldots,g_{\lfloor\frac{d}{2}+1\rfloor})$ (writen as row vectors) and the $f$-vectors $(f_{-1},\ldots, f_{d-1})$ of boundaries of simplicial polytopes of dimension $d$. }

\medskip

In the literature we can find certain results that can be regarded as corollaries of the expresion of McMullen correspondence in Equation \eqref{eq.McM}. We will discuss briefly two examples in the following.

According to \cite{B1}, T. S. Motzkin in the 1950's and D. Welsh in the 1970's seem to have proposed, independently, the following:

\medskip

\noindent \textbf{Unimodality Conjecture.} \emph{The $f$-vector of certain simplicial complexes is always an unimodal sequence. }

\medskip

\noindent A. Bj\"orner disproved this conjecture in \cite{B1} using a counterexample of dimension 24 with $2.6\times 10^{11}$ vertices and, after that,  he proved that such sequences are \emph{almost unimodal} in \cite{B2}. This author also conjectured that the matrix $\bm M_d$ is totally non-negative. This fact was later proved by M. Bj\"orklund and A. Engstr\"om in \cite{BE} and, in a different way, by S. R. Gal in \cite{G2}.

In a similar fashion, we can also find the \emph{upper and lower bound theorems} (we refer the reader \cite[Corollary 8.38]{Z}).

Let us note the following:

\begin{rem}[McMullen correspondence in terms of Riordan matrices] \label{rem.RMcM} Taking our definition of the $g$-vector (the complete vector, not its \emph{first half}), we can state an analogue of McMullen correspondence, given by
$$\bm g\mapsto (T((1-x)^{n+2}\mid 1-x))^{-1}(T(1-x\mid 1))^{-1}\bm g. $$

Doing so, we could take advantage of the fact that some properties for Riordan matrices, such as non-negativity and unimodality of the images, have been explored in the literature (see, for example, \cite{CLW, CWZ, Zhu}).

\end{rem}

Now, let us discuss some aspects about positivity of the $h$-vectors, $g$-vectors and $\gamma$-vectors in connection to Riordan matrices. The $h$-vector of a simplicial complex, in some important cases, has combinatorial meaning (that is, its entries are actually \emph{counting something}), as it is the case of partitionable simplicial complex or of Macaulay complexes  (we refer the reader to \cite{Z}). So, there are some papers in the literature studying classes of simplicial complexes with a positive $h$-vector.  As explained in the first paragraph of the abstract of \cite{BCPT}, for some symmetric $h$-vector with palindromic coefficients, unimodality is equivalent to corresponding to a nonnegative g-vector. Also, a sufficient but not necessary condition for unimodality is having a nonnegative $\gamma$-vector

\begin{rem}[Positivity of $h$-vectors, $g$-vectors and $\gamma$-vectors in terms of Riordan matrices] \label{rem.pos} The  Kruskal-Katona Theorem \cite[Theorem 8.32]{Z} characterizes the set of sequences corresponding to the $f$-vector of a simplicial complex. One could study the problem of finding the vectors $\bm f=(f_{-1},\ldots,f_d)^T$ satisfying the conditions in the Kruskal-Katona Theorem and such that the $h$-vector, the $g$-vector or the $\gamma$-vector is positive. To do this, we can use the Riordan correspondence appearing in  Theorem \ref{theo.changes}, together with the papers concerning positivity of Riordan matrices available in the literature (see, for example, \cite{CLW}).

\end{rem}







\section{$q$-cone of a simplicial complex} \label{sect.properties}

Given two simplicial complexes (either geometric or abstract) $\mathcal{F}$ and $\mathcal{K}$, there is an operation called the \textbf{join} that produces a third simplicial complex, denoted by $\mathcal{F} * \mathcal{K}$. Let us briefly recall the definition of the join operation,
$$\mathcal{F} * \mathcal{K}=\{ \sigma\cup \tau| \sigma\in \mathcal{F}\cup \{\emptyset\}\text{ and }\tau \in \mathcal{K}\cup \{\emptyset\}\},$$

\noindent and let us refer the reader to \cite[Chapter 0]{H}.

\begin{defi} \label{defi.qcone} For $q\geq 1$, let us denote by $[q]$ the 0-dimensional (either abstract or geometric) simplicial complex consisting of $q$ different elements of dimension 0. We call $\mathcal F*[q]$ the \textbf{$q$-cone} of $\mathcal F$.

\end{defi}

In the case $q=1$, the $1$-cone of $\mathcal{F}$ is usually called the \textbf{cone} of $\mathcal{F}$, while in the case $q=2$, $\mathcal{F}*[2]$ is referred to as the \textbf{suspension} of $\mathcal{F}$ in the literature. The aim of this section is to compute and study some combinatorial and topological properties of this family of spaces


\subsection{$f$-vector and homology groups of the $q$-cone of a simplicial complex} \label{subsect.tf}

We can compute easily the $f$-polynomial of the $q$-cone of any simplicial complex $\mathcal F$ in terms of the $f$-polynomial of $\mathcal F$ and the number $q$:


\begin{lemm}\label{L:f-q-cones} Let $\mathcal{F}$ be any $d$-dimensional simplicial
complex with extended $f$-vector $(f_{-1},f_0,f_1,\cdots,f_d)$, $f$-polynomial $P_{\mathcal{F}}(x)$ and exteded $f$-polynomial $\widetilde P_{\mathcal{F}}(x)$. For any positive integer $q$, the extended $f$-vector of $\mathcal F*[q]$ is 
$$(f_{-1},f_0+q,qf_0+f_1,qf_1+f_2,\ldots,qf_{d-1}+f_d),$$

\noindent its $f$-polynomial is
$$P_{\mathcal{F}*[q]}(x)=(qx+1)P_{\mathcal{F}}(x)+q,$$

\noindent and its extended $f$-polynomial is
$$\widetilde P_{\mathcal{F}*[q]}(x)=(qx+1)\widetilde P_{\mathcal{F}}(x).$$

\end{lemm}


On the other hand, we have the following:

\begin{lemm}\label{thm_homologia_conos} Let $\mathcal{F}$ be any simplicial complex, then $$\tilde{H}_{n}(|\mathcal{F}*[q]|)\simeq \bigoplus_{(q-1)-times}\tilde{H}_{n-1}(|\mathcal{F}|).$$

\noindent where $\widetilde H_n$ denotes the corresponding reduced homology group.

\end{lemm}
\begin{proo}
The join of two finite simplicial complexes $\mathcal F$ and $\mathcal F'$ is homotopy equivalent to the suspension of the smash product of $\mathcal F$ and $\mathcal F'$, given by
$\mathcal F\wedge \mathcal F':= (\mathcal F\times \mathcal F') / (\mathcal F\vee \mathcal F')$ (see \cite[Chapter 0, Exercise 23]{H}). Hence, their reduced homology groups are isomorphic. Moreover, we also have that $\tilde{H}_n(|\mathcal F*[2]|)\simeq \tilde{H}_{n-1}(|\mathcal F|)$ for any simplicial complex $\mathcal F$ (see \cite[Chapter 2, Section 2.1, Exercise 20]{H}).



From this, we deduce $$\tilde{H}_{n}(|\mathcal{F}*[q]||)\simeq \tilde{H}_{n-1}(\mathcal{F}\wedge[q])\simeq \tilde{H}_{n-1}((\mathcal{F}\times [q])/(\mathcal{F}\vee [q]))\simeq \bigoplus_{(q-1)-times}\tilde{H}_{n-1}(|\mathcal{F}|).$$

\end{proo}

\begin{rem} \label{rem.homotopy} We have decided to include the previous result for various reasons; the proof is short, it contains what we need for the following section and it ressembles the Riordan pattern, as we discuss below. Nevertheless, we could have replaced it by a stronger statement: \emph{For every simplicial complex $\mathcal F$, $\mathcal F*[q]$ is homotopically equivalent to the wedge of $(q-1)$ copies of the suspension of $\mathcal F$ ($\mathcal F*[2]$).} This result can be proved directly or using the previous one and the fact that every shellable simplicial complex (see below) is homotopy equivalent to a wedge sum of spheres.  We omit details.

\end{rem}




\subsection{Partitionability}



Let us recall the following definition from \cite{Z}:
\begin{defi} \label{defi.partitionable}Let $\mathcal{F}$ be a pure simplicial complex of dimension $d$ and let us denote by \textnormal{facets}$(\mathcal F)$ to the set of facets in $\mathcal F$. An exact facing $\varphi:\textnormal{facets}(\mathcal{F})\rightarrow \mathcal{F}$ is a map satisfying that for every $\sigma\in \mathcal{F}$, there exists a unique facet $\tau \in \textnormal{facets}(\mathcal{F})$ such that $\varphi(\tau)\subseteq \sigma\subseteq \tau$. If there exists an exact facing for $\mathcal{F}$, it is said that $\mathcal{F}$ is a \textbf{partitionable} simplicial complex.
\end{defi}

\begin{theo}\label{thm_join_particionable} Let $\mathcal{F}$ be a finite pure simplicial complex that is partitionable and let $q$ be a positive integer number. Then $\mathcal{F} * [q]$ is a partitionable simplicial complex.
\end{theo} 
\begin{proo}
Let $\varphi:\text{facets}(\mathcal{F})\rightarrow \mathcal{F}$ be an exact facing. Set $\mathcal{F}'=\mathcal{F}*[q]$ and define $\varphi':\text{facets}(\mathcal{F}')\rightarrow \mathcal{F}' $ by $\varphi'(\sigma\cup i)=\varphi(\sigma)$ if $i=1$ and $\varphi'(\sigma\cup i)=\varphi(\sigma)\cup i$ if $1<i\leq q$, where $\sigma\in \text{facets}(\mathcal{F})$. Let us prove that $\varphi'$ is an exact facing by studying cases. 

Suppose that $\tau\in \mathcal{F}\subset \mathcal{F}'$. By hypothesis, there exists a unique facet $\sigma\in \text{facets}(\mathcal F)$ such that $\varphi(\sigma)\subseteq \tau \subseteq \sigma$. It is clear that $\tau \subseteq \sigma \cup i$ for every $1\leq i\leq q$. On the other hand, by construction we get that $\varphi'(\sigma \cup i)\subseteq \tau$ just for $i=1$ because $\varphi'(\sigma \cup 1)=\varphi(\sigma)$ while $\varphi'(\sigma \cup i)=\varphi(\sigma)\cup i \nsubseteq \tau$ for every $1<i\leq q$. If there exists $\sigma'\cup 1$ in $\textnormal{facets}(\mathcal F')$ such that $\sigma\neq \sigma'$ with $\varphi'(\sigma'\cup 1)\subseteq \tau \subseteq \sigma'\cup 1$, we get a contradiction with the uniqueness of $\sigma$. Thus, for every $\tau \in \mathcal{F}\subset \mathcal{F}'$ there exists a unique facet $\sigma\cup i$ in $\text{facets}(\mathcal{F}')$ such that $\varphi'(\sigma\cup 1)\subseteq \tau \subseteq \sigma\cup 1$.

Suppose now that $\tau \in \mathcal{F}'\setminus \mathcal{F}$, so $\tau=\gamma\cup i$ for some $\gamma\in \mathcal{F}$ and $1\leq i\leq q$.  It is clear that $\gamma\cup i\subseteq \sigma \cup j$ if and only if $j=i$ and $\gamma\subseteq \sigma$ (where $\sigma\in \text{facets}(\mathcal{F})$). By hypothesis, there exists  a unique facet $\sigma\in \text{facets}(\mathcal{F})$ such that $\varphi(\sigma)\subseteq \gamma\subseteq \sigma$. Hence, $\varphi'(\sigma\cup i)\subseteq \gamma\cup i\subseteq \sigma\cup i$. If there exists $\sigma\neq \sigma'\in \text{facets}(\mathcal{F})$ such that $\varphi'(\sigma'\cup i)\subseteq \gamma\cup i\subseteq \sigma'\cup i$, then $\varphi(\sigma')\subseteq \gamma \subseteq \sigma'$, which entails a contradiction. Thus, $\sigma\cup i$ is the unique facet in $\text{facets}(\mathcal{F}')$ satisfying $\varphi'(\sigma\cup i)\subseteq \gamma\cup i \subseteq \sigma\cup i$. 
\end{proo}

One interesting property of partitionable simplicial complexes 
is that the exact facings have a combinatorial meaning. In particular,  if $\varphi$ is an exact facing of a $d$-dimensional simplicial complex $\mathcal{F}$, then its $h$-vector is given by 
$$h_j=|\{ \sigma\in \mathcal{F}| \  \sigma\  \textnormal{is a facet and }  \textnormal{dim} (\varphi(\sigma))=j\}|, \ \ j=0,...,d. $$

\noindent For example, see \cite[Chapter III, Proposition 2.3]{Stanley}. As a consequence, the entries in the $h$-vector are all of them positive and we have the following:

\begin{theo}\label{thm_exact_facing_h_vector_join} Let $\mathcal{F}$ be a finite pure $d$-simplicial complex and let $q$ be a positive integer number. Assume that the $h$-vector of $\mathcal{F}$ is $(h_0,\ldots, h_{d})$ and its $h$-polynomial is $Q_{\mathcal F}(x)$. If $\varphi:\textnormal{facets}(\mathcal{F})\rightarrow \mathcal{F}$ is an exact facing for $\mathcal{F}$, then the $h$-vector $(h_0',\ldots, h_{d+1}')$ of $\mathcal{F}*[q]$ is given by 
$$h_j'=h_j+(q-1)h_{j-1}, \ \ j=0,...,d+1 $$

\noindent and its $h$-polynomial $Q_{\mathcal F*[q]}(x)$ is $Q_{\mathcal F}(x)+(q-1)xQ_{\mathcal F}(x). $

\end{theo}

\begin{proo}
It is a direct consequence of the proof of Theorem \ref{thm_join_particionable}, where it is constructed an explicit exact facing for $\mathcal{F}*[q]$ from $\varphi:\text{facets}(\mathcal{F})\rightarrow \mathcal{F}$. 

\end{proo}





\begin{rem}\label{rem:interpretacion.exact.facing.particiones}
For a partitionable pure $d$-dimensional simplicial complex $\mathcal{F}$, the entry $h_i$ of the $h$-vector counts the intervals of height $d-i$ in the corresponding partitioning of the face poset of $\mathcal{F}$, i.e., $\mathcal{X}(\mathcal{F})$, into intervals. Specifically, the intervals are defined as $[\varphi(\sigma),\sigma]$, where $\varphi:\textnormal{facets}(\mathcal{F})\rightarrow \mathcal{F}$ is an exact facing. 
\end{rem}

\subsection{Vertex-decomposability and shellability}

Let us start by recalling the following two concepts:
\begin{defi}\label{df:vertex_deco} Let $\mathcal F$ be a simplicial complex. For every face $\sigma \in \mathcal{F}$, the \textbf{link} and the \textbf{deletion} of $\sigma$ are defined, respectively as:
$$\textnormal{lk}_\mathcal{F}(\sigma)=\{ \tau\in \mathcal{F} |\sigma\cap \tau= \emptyset\ \ \text{and} \ \ \sigma\cup \tau\in \mathcal{F} \}, \qquad\textnormal{del}_\mathcal{F}(x)=\{ \tau\in \mathcal{F}|\sigma\cap \tau= \emptyset\  \}.$$

Suppose that $V=\{x_1,...,x_n\}$ is the vertex set of $\mathcal F$. Then $\mathcal{F}$ is \textbf{vertex-decomposable} if either:
\begin{enumerate}
\item The only facet of $\mathcal{F}$ is $\{x_1,...,x_n \}$, i.e., $\mathcal{F}$ is a simplex, or $\mathcal{F}=\emptyset$.
\item There exists $x\in V$ such that $\textnormal{del}_\mathcal{F}(x)$ and $\textnormal{lk}_{\mathcal{F}}(x)$ are vertex decomposable, and such that every facet of  $\textnormal{del}_\mathcal{F}(x)$ is a facet of $\mathcal{F}$.
\end{enumerate}
\end{defi}

\begin{defi} Let $\mathcal{F}$ be a pure simplicial complex. We say that $\mathcal{F}$ is \textbf{shellable} if its facets can be ordered $F_1,F_2,...,F_s$ such that the following condition holds. Write $\mathcal{F}_j$ for the subcomplex of $\Delta$ generated by $F_1,...,F_j$. Then we require that for all $1\leq i\leq s$, the set of faces of $\mathcal{F}_i$ which do not belong  to $\mathcal{F}_{i-1}$ has a unique minimal element (with respect to inclusion). (When $i=1$, we have $\mathcal{F}_0=\emptyset$, so $\mathcal{F}_1\setminus \mathcal{F}_0$ has the unique minimal element $\emptyset$.)
\end{defi}

As stated in \cite[Chapter III, Section 2]{Stanley}, if a simplicial complex $\mathcal F$ is shellable, then it is partitionable. However, the converse does not hold (see \cite[Example 2.1]{KO}). Additionally, if $\mathcal{F}$ is a pure vertex decomposable simplicial complex, then $\mathcal{F}$ is shellable (see \cite[Corollary 2.9]{Prov80}. Again, the converse fails to be true (\cite[Appendix F]{Sim1994}). We can sum up this situation in the following diagram:
\[\begin{tikzcd}
\textnormal{Vertex-decomposable} \arrow[r, Rightarrow, shift left=1 ex] & \textnormal{Shellable} \arrow[r, Rightarrow,  shift left=1 ex] \arrow[l,Rightarrow, "/" marking, shift left=1 ex]&  \textnormal{Partitionable} \arrow[l,Rightarrow, "/" marking, shift left=1 ex]
\end{tikzcd}
\]

In light of the previous discussion, it suffices to show that a simplicial complex is vertex decomposable to prove that it is also partitionable. However, we have chosen to explicitly demonstrate that $\mathcal{F}*[q]$ is partitionable, as we believe this approach is illustrative and allows us to determine the corresponding $h$-vector. Nonetheless, the following result is stronger, and Theorem \ref{thm_join_particionable} can be viewed as a direct corollary of Theorem \ref{thm_vertex_decomposable_join}.

\begin{theo}\label{thm_vertex_decomposable_join} Suppose that $\mathcal{F}$ is a vertex-decomposable finite simplicial complex. Then $\mathcal{F} * [q]$ is a vertex-decomposable finite simplicial complex.
\end{theo}
\begin{proo}
We argue by induction on $q$. Suppose $q=1$ and set $\mathcal{F}_1=\mathcal{F}*[1]$. Suppose that $\mathcal{F}$ is a simplex, then $\mathcal{F}_1$ is a simplex and therefore it is vertex decomposable. Suppose that there exists a vertex $x$ such that $\text{del}_\mathcal{F}(x)$ and $\text{lk}_{\mathcal{F}}(x)$ are vertex decomposable, and such that every facet of $\text{del}_\mathcal{F}(x)$ is a facet of $\mathcal{F}$. (*) Note that $\text{lk}_{\mathcal{F}_1}(x)=\text{lk}_{\mathcal{F}}(x)*\{1\}$ is isomorphic to $\mathcal{E}=\text{lk}_{\mathcal{F}}(x)*x$. By hypothesis, $\text{lk}_{\mathcal{E}}(x)=\text{lk}_\mathcal{F}(x)$ is vertex decomposable. Moreover, $\text{del}_{\mathcal{E}}(x)=\text{lk}_\mathcal{F}(x)$ so it is vertex decomposable. This gives that $\text{lk}_{\mathcal{F}_1}(x)$ is vertex decomposable. On the other hand, $\text{del}_{\mathcal{F}_1}(x)=\text{del}_{\mathcal{F}}(x)*\{ 1\}$, where $\text{del}_{\mathcal{F}}(x)$ is vertex decomposable and trivially every facet of $\text{del}_{\mathcal{F}_1}(x)$ is a facet of $\mathcal{F}_1$. There are two options. (1) If $\text{del}_{\mathcal{F}}(x)$ is a simplex or the empty set, then $\text{del}_{\mathcal{F}_1}(x)$ is vertex decomposable and every facet of $\text{del}_{\mathcal{F}_1}(x)$ is a facet of $\mathcal{F}_1$. (2) If $\text{del}_{\mathcal{F}}(x)$ is not a simplex or the empty set, then there exists $y\in \text{del}_{\mathcal{F}}(x)$ such that $\text{lk}_{\text{del}_{\mathcal{F}}(x)}(y)$ and $\text{del}_{\text{del}_{\mathcal{F}}(x)}(y)$ are vertex decomposable and every facet of $\text{del}_{\text{del}_{\mathcal{F}}(x)}(y)$ is a facet of $\text{del}_{\mathcal{F}}(x)$. Note that $\text{lk}_{\text{del}_{\mathcal{F}_1}(x)}(y)=\text{lk}_{\text{del}_{\mathcal{F}(x)}(x)}(y)*\{ 1\}$ is isomorphic to $\text{lk}_{\text{del}_{\mathcal{F}}(x)}(y)* y$. Hence, we only need to repeat the same arguments used in $(*)$ to get that $\text{lk}_{\text{del}_{\mathcal{F}}(x)*\{ 1\}}(y)$ is vertex decomposable. It remains to check whether $\text{del}_{\text{del}_{\mathcal{F}}(x)*\{ 1\}}(y)=\text{del}_{\text{del}_{\mathcal{F}}(x)}(y)*\{ 1\}$ is vertex decomposable. As before, there are two options and we just need to repeat the same arguments. Since $\mathcal{F}$ is a finite set, this process ends with an empty set or a simplex.

Suppose the result holds for $q$ and consider $q+1$. Set $\mathcal F_{q}=\mathcal{F}*[q]$ and $\mathcal 
 F_{q+1}=\mathcal{F}*[q+1]$. Consider $x=\{q+1 \}$. Then $\text{lk}_{\mathcal{F}_{q+1}}(x)=\mathcal{F}$ and $\text{del}_{\mathcal F_{q+1}}(x)=\mathcal F_q$. Both simplicial complexes are vertex decomposable and every facet of $\text{del}_{\mathcal F_{q+1}}(x)$ is a facet of $\mathcal F_{q+1}$, which gives that $\mathcal F_{q+1}$ is vertex decomposable.
\end{proo}

\section{Iterated $q$-cones}\label{sect.qc}

Let us recall that, for every positive integer, we denote by $[k]$ to the 0-dimensional simplicial complex consisting of $k$ vertices.

\begin{defi} \label{defi.qsec} For every positive integers $m,q$, we define the sequence $(\Delta_{m,q}^{(n)})_{n\in\mathbb N}$ by
$$\Delta_{m,q}^{(0)}=[m], \qquad \forall n\geq 1,\; \Delta_{m,q}^{(n)}=\Delta_{m,q}^{(n-1)}*[q].  $$

\noindent We refer to $\Delta_{m,q}^{(n)}$ as the \textbf{$n$-th $q$-cone with initial condition $m$}.
\end{defi}

\begin{figure}[h!]
\centering
\includegraphics[width=0.7\textwidth]{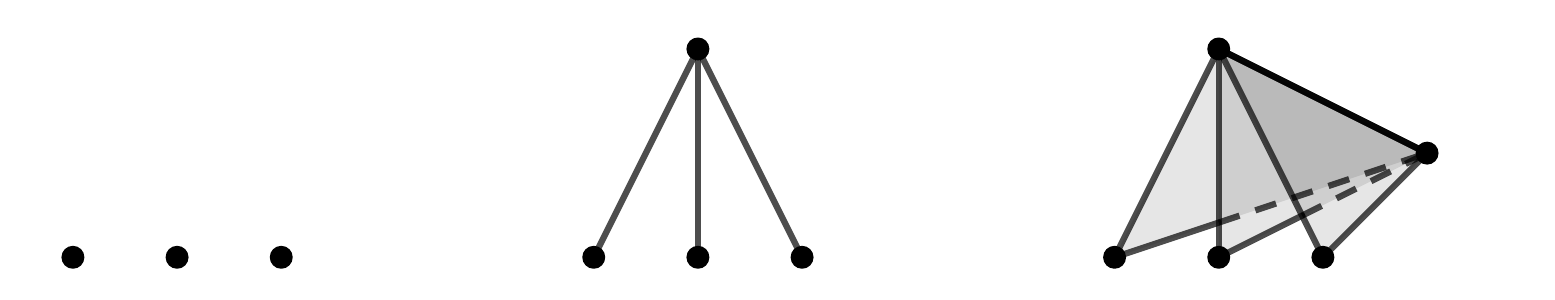}
\caption{Geometric realizations of the simplicial complexes $\Delta_{31}^{(0)}$, $\Delta_{31}^{(1)}$ and $\Delta_{31}^{(2)}$.}

\end{figure}


\subsection{$f$-vector and $h$-vector of the $n$-th $q$-cone}\label{sub:f-vector_matrix_h-vector}



In this subsection, we obtain Riordan matrices from the $f$-vectors and $h$-vectors of $\Delta_{m,q}^{(n)}$.

\begin{theo} \label{theo.muchos} Let $\bm F_{m,q}=(f_{n,k})_{n,k\geq0}$ be the infinite matrix such that the
$n$-row $(f_{n0},\ldots, f_{nn},0,\ldots)$ of $\bm F_{m,q}$ is just the $f$-vector of $\Delta^{(n)}_{m,q}$ or, equivalently, such that the generating function of the $n$-row is the $f$-polynomial of  $\Delta^{(n)}_{m,q}$.

Let us denote by $P_0(x),P_1(x),P_2(x),\ldots$ to the $f$-polynomials of $\Delta^{(0)}_{m,q},\Delta^{(1)}_{m,q},\Delta^{(2)}_{m,q},\ldots$ respectively (and thus, the generating functions of the rows in $\bm F_{m,q}$).

Then

\medskip

\noindent (a)  The sequence of polynomials $P_0(x),P_1(x),P_2(x),\ldots$ is determined by the recurrence
$$P_0(x)=m,\qquad P_n(x)=(qx+1)P_{n-1}(x)+q.$$

\noindent (b) $\bm F_{m,q}$ is a Riordan matrix. Specifically,
\[
\bm F_{m,q}=T\left(\frac{m+(q-m)x}{q(1-x)}\Bigm|\frac{1-x}{q}\right).
\]

\noindent (c) The entries in $\bm F_{m,q}$ can be described iteratively by
$$\begin{array}{l l} \text{(description of the first column)} &f_{0,0}=m,\qquad \forall n\geq1, \; f_{n,0}=f_{n-1,0}+q, \\
\text{(A-sequence rule)}& \forall n,k\geq1,\; f_{n,k}=f_{n-1,k}+qf_{n-1,k-1}.\end{array}$$

\medskip

\noindent (d) $\forall\ n,k\geq 0$,
\[
f_{n,k}=\left[q\binom{n}{k+1}+m\binom{n}{k}\right]q^k.
\]

\noindent (e) The solution of the recurrence in Statement (a) is
$$P_n(x)=m(qx+1)^{n+1}+q\sum_{j=0}^n(qx+1)^j.$$


\end{theo}

\begin{proo} Statement (a) is a direct consequence of the fact that $\Delta^{(0)}_{m,q}=[m]$ and of  Lemma \ref{L:f-q-cones}.

To prove Statement (b), using now Theorem 5 in \cite{poly} we obtain that $\bm F_{m,q}$ is the
Riordan matrix, which appears to \emph{the right side of the
vertical line} and \emph{below the horizontal line} in
\begin{small}
\begin{equation}\label{eq.line}
\left(
  \begin{array}{c|cccccc}
\frac{m}{q}&&&&&&\\ \hline
1&m&&&&&\\
1&m+q&mq&&&&\\
1&m+2q&2mq+q^2&mq^2&&&\\
1&m+3q&3mq+3q^2&3mq^2+q^3&mq^3&&\\
1&m+4q&4mq+6q^2&6mq^2+4q^3&4mq^3+q^4&mq^4&\\
\vdots&\vdots&\vdots&\vdots&\vdots&\vdots&\ddots\\
  \end{array}
\right),
\end{equation}
\end{small}
which corresponds to  $T(\alpha\mid \omega)$ for
\[
\alpha(x)=\frac{m}{q}+\frac{x}{1-x}, \qquad \omega(x)=\frac{1}{q}-\frac{1}{q}x,
\]

To prove Statement (c) note that, by definition of $\Delta^{(n)}_{m,q}$, for ever non-negative integer $n$,  $\Delta^{(n)}_{m,q}$ has $m+nq$ vertices. After this, we can use Theorem 11 in \cite{teo} that directly induces the recurrence in Statement (c).

To prove Statement (d), we  use the product in the Riordan group to obtain the following identity:
\[
\bm F_{m,q}=T\left(\frac{m+(q-m)x}{q(1-x)}\Bigm|1\right)T\left(1\Bigm|1-x\right)T\left(1\Bigm|\frac{1}{q}\right).
\]

Note that  $T(1\mid 1-x)$ is the Pascal triangle. The polynomials $S_n(x)$ corresponding to the generating functions of the rows of the Pascal triangle (in a similar fashion to the polynomials $P_n$ corresponding to the generating functions of the rows of $\bm F_{m,q}$)
are just $S_n(x)=(1+x)^n$. Using Propositions 12 and 18 in
\cite{poly} we get that
\[
P_n(x)=q\sum_{k=0}^{n-1}S_k(qx)+mS_n(qx)=\sum_{k=0}^n\left[q\binom{n}{k+1}+m\binom{n}{k}\right]q^kx^k
\]
and so we obtain (d).

Statement (e) is a consequence of Statement (d) or, equivalently, can be proved straightforward by induction from Statement (a).


\end{proo}

To illustrate the previous theorem, the reader may find the matrix $\bm F_{1,1}$ in Equation \eqref{eq.F11} and the matrix $\bm F_{2,2}$ in Equation \eqref{eq.F22} (Subsection \ref{subsect.exam}). Also, in each case ($m=q=1$ and $m=q=2$), the explicit formulas for the entries of the $f$-vectors (as they appear in Statement (d)), can be found in Equations \eqref{eq.11fnk} and  \eqref{eq.22fnk}.

Recall that the extended $f$-vector of $\mathcal{F}$ is just $(f_{-1},f_0,f_1,\cdots,f_d)$ where $(f_0,\ldots, f_d)$ is the (standard) $f$-vector and $f_{-1}=1$. In a similar fashion, if $P_{\mathcal F}(x)$ is the corresponding $f$-polynomial, then the extended $f$-polynomial is obtained as $\widetilde P_{\mathcal F}=1+xP_{\mathcal F}(x)$. Using arguments as in \cite{2ways} and the previous result, we obtain:




\begin{coro} \label{coro.extendedF}  Let us define the matrix $\widetilde{\bm F}_{m,q}=(d_{nk})$ such that $f_{00}=m/q$ and such that, for all $n\geq 1$, $(d_{n0},\ldots, d_{nn})$ is the extended $f$-vector of $\Delta^{(n-1)}_{m,q}$ (note that $d_{nk}$ is the number of $k-1$ dimensional faces of $\Delta^{(n-1)}_{m,q}$, in this setting). Then
\[
\widetilde{\bm F}_{m,q}=T\left(\frac{m+(q-m)x}{q^2}\Bigm|\frac{1-x}{q}\right).
\]


\end{coro}

Note that
\begin{small}
\[
\widetilde{\bm F}_{m,q}= \left(
  \begin{array}{ccccccc}
\frac{m}{q}&&&&&&\\
1&m&&&&&\\
1&m+q&mq&&&&\\
1&m+2q&2mq+q^2&mq^2&&&\\
1&m+3q&3mq+3q^2&3mq^2+q^3&mq^3&&\\
1&m+4q&4mq+6q^2&6mq^2+4q^3&4mq^3+q^4&mq^4&\\
\vdots&\vdots&\vdots&\vdots&\vdots&\vdots&\ddots\\
  \end{array}
\right)
\]
\end{small}

\noindent is obtained  form the matrix in  Equation \eqref{eq.line} by deleting the vertical and horizontal lines. That is, $\widetilde{\bm F}_{m,q}$ is obtained from $\bm F_{m,q}$ by properly adding one column on the left and one row on the top or, equivalently, $\bm F_{m,q}$ is obtained from $\widetilde{\bm F}_{m,q}$ by deleting the 0-row and the 0-column.





The expression in Statement (c) in the theorem above is not the only possibility of interest for describing the number of $k$-dimensional faces of the simplicial complex $\Delta_{m,q}^{(n)}$ as a linear combination of the number of faces of the previous elements in the sequence $\Delta_{m,q}^{(n-s)}$, for $s=1, \ldots, n$. For example, we have

\begin{theo}
Given a nonnegative integer $s$, 
\begin{equation}\label{e:A-seq}
f_{n,k}=\sum_{j=0}^{n-k}q^{s-j}\binom{s}{j}f_{n-s, k-s+j},
\end{equation}  
and
\begin{equation}\label{e:g-seq}
f_{n,k}=\sum_{j=0}^{n-k}q^{k}\binom{s}{j}f_{n-s-j, k-s}.
\end{equation} 
\end{theo}

\begin{proo}
Since the $A$-sequence of $\bm F_{m,q}$ is $A(x)=q+x$ and $\omega(x)=\frac{1-x}{q}$, using the notation in Theorem 1 in \cite{double} (or Theorem 3.3 in \cite{iden}) we directly get the result. 
\end{proo}

The equality (\ref{e:A-seq}) describes the number of $k$-dimensional faces of the simplicial complex $\Delta_{m,q}^{(n)}$ as a linear combination of different numbers corresponding to the faces of various dimensions of the simplicial complexes $\Delta_{m,q}^{(n-s)}$. However, in the equality (\ref{e:g-seq}), the same number is expressed as a linear combination of the number of $(k-s)$-dimensional faces of different iterations of our $q$-cone construction with the initial condition $m$.


\begin{theo} \label{theo.hmq}Let us denote by $\bm H_{m,q}=(h_{ij})_{0\leq i,j<\infty}$ such that:

\begin{itemize}

\item the 0-row is determined by imposing that $h_{00}=(m-1)/(q-1)$ and, for $j\geq 1$, $h_{0j}=0$,

\item  the $i$-row, for $i\geq 1$, is determined by imposing that $h_{i0}+h_{i1}x+\ldots+h_{ii}x^i$ is the $h$-polynomial of $\Delta_{m,q}^{(i-1)}$ and, for $j\geq i+1$, $h_{ij}
=0$.

\end{itemize}

\noindent  Then, for all $m,q\geq 2$, $\bm H_{m,q}=\widetilde{\bm F}_{m-1,q-1}$ (see Corollary \ref{coro.extendedF}, for the definition of $\widetilde{\bm F}_{m-1,q-1}$).

\end{theo}

\begin{proo} This result is a consequence of Lemma \ref{L:f-q-cones} and Theorem \ref{thm_exact_facing_h_vector_join}.

\end{proo}

\begin{rem} As a consequence of the previous theorem, note that the $h$-vectors of these spaces are, in fact, $f$-vectors. This means that they satisfy, for instance, the Kruskal-Katona conditions. These kinds of questions have been of interest to many authors (see \cite{BiermanAdam,BFS} and the references therein).
\end{rem}
To conclude this section, we establish a matrix relation between $\bm{H}_{m,q}$ and $\widetilde{\bm{F}}_{m,q}$, leveraging the group properties of Riordan matrices.

\begin{prop} Let $\bm H_{m,q}$ and $\widetilde{\bm F}_{m,q}$ be the matrices described before. Then,
$$\bm H_{m,q}=\underbrace{T\left(\frac{m-1+(q-m)x}{m+(q-1-m)x}\frac{q}{q-1}\frac{q-x}{q-1}\Bigm| \frac{q-x}{q-1}\right)}_{\bm D_L}\widetilde{\bm F}_{m,q}, $$

$$\bm H_{m,q}=\widetilde{\bm F}_{m,q}\underbrace{T\left(\frac{q(m-1)+(q-1)x}{(q-1)(m+x)}\frac{q}{q-1}\Bigm|\frac{q}{q-1}\right)}_{\bm D_R}. $$

In the regular case, it is $m=q$,
$$\bm D_L=T\left(\frac{q}{q-1}\Bigm|\frac{q-x}{q-1}\right),\qquad \bm D_R=T\left(\frac{q}{q-1}\Bigm|\frac{q}{q-1}\right).$$

Moreover, since $\widetilde{\bm F}_{11}$ is the Pascal Triangle, we have
$$\widetilde{\bm F}_{q,q}=T\left(\frac{1}{q}\Bigm|\frac{1+(q-1)x}{q}\right) \widetilde{\bm F}_{1,1},\qquad \widetilde{\bm F}_{q,q}= \widetilde{\bm F}_{1,1}T\left(\frac{1}{q}\Bigm|\frac{1}{q}\right).$$

\end{prop}

\subsection{The Euler characteristic and other arithmetical invariants}
Let us introduce the Euler characteristic and some of its properties:
\begin{defi}[Euler characteristic as a combinatorial invariant] \label{defi.ec} Given a $d$-dimensional simplicial complex $\mathcal{F}$ with f-vector $(f_{0},\ldots,f_d)$, the
\textbf{Euler characteristic} of $\mathcal F$ is defined as
\[
\chi(\mathcal{F})=\sum_{k=0}^{d}(-1)^kf_k.
\]
\end{defi}

Consider the geometric realization $|\mathcal{F}|$ of $\mathcal{F}$. It is well known that $\chi(\mathcal{F})$ depends on $|\mathcal{F}|$. This is due to the fact that $\chi(\mathcal{F}) = \sum_{i=0} (-1)^i \beta_i$, where $\beta_i = \textnormal{rank} \ H_i(|\mathcal{F}|)$ (here we are considering singular homology with coefficients in a field). Specifically, this result indicates that the Euler characteristic is a homotopy invariant (see \cite[Theorem 2.44]{H}).


\begin{rem} \label{rem.pr} Note that $\chi(\mathcal{F})$ can be expressed as a matricial
product
\begin{small}
\[
\chi(\mathcal{F})=(f_0,f_1,\cdots,f_d,0,\cdots)\begin{smallmatrix} \left(%
\begin{array}{c}
  1 \\
  -1 \\
  1\\
  -1 \\
  \vdots \\
   (-1)^n \\
    \vdots \\
\end{array}%
\right)
\end{smallmatrix}
\]
\end{small}

\noindent and the column on the right does not depend on the complex $\mathcal{F}$.

\end{rem}

This observation and basic properties of Riordan arrays allow us
to get at once the following:

\begin{theo} \label{theo.vectorChi}
Let $m,q\geq1$ be positive integers and let
$\chi_{m,q}^{(n)}$ denote $\chi(\Delta_{m,q}^{(n)})$. Consider the sequence
$\left(\chi_{m,q}^{(n)}\right)_{n\in\N}$ and denote by
$\chi_{m,q}(x)=\sum_{n\geq0}\chi_{m,q}^{(n)}x^n$ the generating
function of $\left(\chi_{m,q}^{(n)}\right)_{n\in\N}$. Then
\[
\chi_{q,m}(x)=\frac{m+(q-m)x}{(1-x)(1-(1-q)x)}
\]
and, consequently,
\[
\chi_{q,m}^{(n)}=1+(m-1)(1-q)^n
\]
\end{theo}
\begin{proo}
Using the above interpretation of the Euler characteristic we have
the following matricial equality
\[
\bm F_{m,q}\left(\frac{1}{1+x}\right)=T\left(\frac{m+(q-m)x}{q(1-x)}\Big|\frac{1-x}{q}\right)\left(\frac{1}{1+x}\right)=
\chi_{q,m}(x),
\]
but
\[
T\left(\frac{m+(q-m)x}{q(1-x)}\Big|\frac{1-x}{q}\right)\left(\frac{1}{1+x}\right)=
\frac{m+(q-m)x}{q(1-x)}\frac{1-x}{1-x+qx}=\frac{m+(q-m)x}{(1-x)(1-(1-q)x)}.
\]
Now expanding into power series this rational function we get
\begin{equation} \label{eq.eulercarqc}
\chi_{q,m}(x)=\sum_{k\geq0}(1+(m-1)(1-q)^n)x^n.
\end{equation}
\end{proo}

It is well known that the Euler characteristic of the simplex of each dimension is 1. This is coherent with the result above because
$$\bm F_{1,1}\left(\frac{1}{1+x}\right)=\frac{1}{1-x}. $$

It is worth noting the uniqueness of the Euler Characteristic. By uniqueness, we mean that if there exists a function assigning a number to each simplicial complex $\mathcal{F}$ based on the number of simplices in $\mathcal{F}$, this function is either the Euler Characteristic or closely related to it (see, for instance, \cite{Brown, levitt}). Specifically, in \cite{Brown}, it is proven that \textit{the only subdivision invariants (linear or not) of $\mathcal{F}$
are the Euler characteristic and the dimension}. Recently, in \cite{LMP}, using the Riordan language, it has been proven that multiples of the Euler Characteristic are the unique invariants under barycentric subdivisions within the class of all $n$-simplices. 

Based on the previous comments, one might wonder what series $\mu=\mu_0+\mu_1+\mu_2x^2+\ldots$ would saltisfy
$$\bm F_{m,q}(\mu)=\frac{1}{1-x} $$

\noindent because this would provide a common linear arithmetic relation for the number of faces of the simplicial complex $\Delta_{m,q}^{(n)}$. To find this series, we only need to solve the following:
$$\bm F_{m,q}(\mu)=\frac{1}{1-x} $$

\noindent or, equivalently,
$$\mu=\bm F_{m,q}^{-1}\left(\frac{1}{1-x}\right)=T\left(\frac{q}{m+x}\mid q+x\right)\left(\frac{1}{1-x}\right)=\frac{\frac{q}{mx}}{q+x}\frac{1}{1-\frac{x}{q+x}}=\frac{1}{m+x}=\frac{1}{m}\left(\frac{1}{1-\left(-\frac{x}{m}\right)}\right), $$

\noindent that is, the coefficients corresponding to $\mu$ are $\frac{1}{m},-\frac{1}{m^2},\frac{1}{m^3},\ldots$ Surprisingly, they depend only on $m$ and not on $q$. We have proved the following:

\begin{coro} \label{coro.meuler} Given $m\geq 1$, for every $q\geq 1$ and every $n\geq 0$, we have
$$\frac{1}{m}f_{n0}-\frac{1}{m^2}f_{n1}+\frac{1}{m^3}f_{n2}-\ldots+(-1)^n\frac{1}{m^{n+1}}f_{nn}=\sum_{k=0}^n(-1)^k\frac{1}{m^{k+1}}f_{nk}=1, $$


\noindent where $f_{nk}$ is the number of $k$-dimensional faces of the simplicial complex $\Delta_{m,q}^{(n)}$.

\end{coro}

In the following, let us denote by $M_m(x)=\frac{1}{m}-\frac{1}{m^2}x+\ldots=\sum_{k=0}^n(-1)^k\frac{1}{m^{k+1}}x^k=\frac{1}{m+x}$. Note that, for $m=1$, we get $M_1(x)=\frac{1}{1+x}$.

Moreover, the previous equality can be written as
$$m^nf_{n0}-m^{n-1}f_{n1}+m^{n-2}f_{n2}-\ldots+(-1)^nf_{nn}=\sum_{k=0}^n(-1)^km^{n-k}f_{nk}=m^{n+1}. $$

\subsection{A topologically generated Riordan matrix} \label{subsect.betti}

Note that one of the ways to interpret the matrix $\bm F_{m,q}$ is the
following: \textit{$f_{n,k}$ is the rank of the free abelian group
$C_k(\Delta^{(n)}_{q,m})$ of simplicial $k$-chains in
$\Delta^{(n)}_{q,m}$}. By using Betti numbers, we also obtain Riordan matrices derived from topological invariants.
\begin{theo} \label{theo.matrixB}
Let $q\geq2$ and let 
$\tilde{\be}_{n,k}=\textnormal{rank}\  \tilde{H}_k(\Delta^{(n)}_{q,m})$. Then the matrix $\bm B_{m,q}=(\tilde{\be}_{n,k})_{n,k\geq0}$ is a Riordan
matrix. Specifically,
\[
\bm B_{m,q}=T\left(\frac{m-1}{q-1}\Bigm|\frac{1}{q-1}\right).
\]
Consequently
\[
\tilde{\be}_{n,k}=\left\{%
\begin{array}{ll}
    (m-1)(q-1)^n, & \hbox{if $n=k$;} \\
    0, & \hbox{if $n\neq k$.} \\
\end{array}%
\right.
\]
If $q=1$, then $\bm B_{m,1}$ is the matrix whose entries are all null except for
$\tilde{\be}_{0,0}=m-1$ which corresponds to the fact that 
$q=1$ and $m\geq1$ imply that $\Delta^{(n)}_{q,m}$ is contractible for
$n\geq1$.
\end{theo}

\begin{proo} Using iteratively Theorem \ref{thm_homologia_conos}, we obtain the recurrence
\[
\tilde{\be}_{n,k}=(q-1)\tilde{\be}_{n-1,k-1} \quad \text{for} \
n,k\geq1, \qquad \tilde{\be}_{0,0}=1\quad \text{and} \quad
\tilde{\be}_{n,0}=0, \quad \text{for} \quad n\geq1
\]
because $|\Delta_{q,m}^{(n)}|$ is path-connected if $n\geq1$ (this is not the case for $n=0$ and $m>1$). All of the above 
implies that if $\bm B_{m,q}=(\tilde{\be}_{n,k})_{n,k\geq0}$, where
$\tilde{\be}_{n,k}=\textnormal{rank}\  \tilde{H}_k(\Delta^{(n)}_{q,m})$,
then
\begin{small}
\[
\bm B_{m,q}=T\left(\frac{m-1}{q-1}\Bigm|\frac{1}{q-1}\right)=
\left(
\begin{array}{ccccc}
  m-1 &  &  &  &  \\
  0 & (q-1)(m-1) &  &  &  \\
  0 & 0 & (q-1)^2(m-1) &  &  \\
  \vdots & \vdots & \vdots & \ddots &  \\
  0 & 0 & 0 & \cdots & (q-1)^n(m-1)) \\
  \vdots & \vdots & \vdots & \vdots & \ddots
\end{array}%
\right)
\]
\end{small}
\end{proo}

\begin{rem} In a similar fashion to Remark \ref{rem.homotopy}, let $m,q$ and $n$ be positive integers. Then $\Delta_{m,q}^{(n)}$ is homotopically equivalent to the wedge sum of $(m-1)(q-1)^n$ copies of $n$-dimensional spheres.
\end{rem}

Notice that, given the fact that the Euler characteristic can be defined as in Definition \ref{defi.ec} or using Betti numbers, we can derive Theorem \ref{theo.vectorChi} from Theorem \ref{theo.matrixB}.

\subsection{Connection to Alexandroff spaces or partially order sets} \label{subsect.qalexandroff}
In this section, we derive some of the previously defined matrices using alternative methods. Let us recall the combinatorial analogue of the join within the framework of the Alexandrofff spaces, that is, the \textbf{non-Hausdorff join}. Given two Alexandroff spaces $X$ and $Y$, the non-Hausdorff join $X\circledast Y$  is the disjoint union $X\sqcup Y$ keeping the given ordering within $X$ and $Y$ and setting $x\leq y$ for every $x\in X$ and $y\in Y$. It is straightforward to verify that $\mathcal{K}(X\circledast Y)=\mathcal{K}(X)*\mathcal{K}(Y)$.  We now focus on the study of the properties of the non-Hausdorff join when iterating infinitely  many times.  
\begin{theo}\label{thm:contractible_general}
Let  $\{X_n \}_{n\in \mathbb{N}}$ be a sequence of non-empty finite spaces and $\Delta=X_{0}\circledast X_{1}\circledast \cdots \circledast X_{i} \circledast X_{i+1}\circledast\cdots$. Then $\Delta$ is contractible.
\end{theo}
\begin{proo}
For each finite space $X_i$, consider a point $x_1^i$ such that $\text{ht}(x_1^i)=0$. Define $c:\Delta\rightarrow \Delta$ by $c(x):=x_1^0$ for every $x\in X$, i.e., $c$ is a constant map. Let $i:\{x_1^0 \}\rightarrow \Delta$ be the inclusion map. It is evident that $c \circ i$ is the identity map $\textnormal{id}_{x_1^0}:\{x_1^0 \}\rightarrow \{x_1^0 \}$. We prove that $i \circ c$ is homotopic to the identity map $\textnormal{id}_{\Delta}:\Delta\rightarrow \Delta$. To do this, consider $S:\Delta\rightarrow \Delta$ defined by $S(y):=x_1^{i+1}$, where $y\in X_i$. Clearly, $S$ is well defined, and we now verify the continuity of $S$ by studying different cases: (i) If $y\leq z$ in $X_i$ for some $i\in \mathbb{N}$, then $S(y)=x^{i+1}_1=S(z)$ (ii) If $y\in X_i$ and $z\in X_j$ such that  $y<z$. Then, we have $i<j$ and, consequently, $S(y)=x_1^{i+1}<x_1^{j+1}=S(z)$. Hence, $S$ is continuous because it is order-preserving. To conclude, let us suppose $y\in \Delta$, which implies that there exists $i\in \mathbb{N}$ such that $y\in X_i$. We obtain that
$$\textnormal{id}_\Delta(y)=y<x_1^{i+1}=S(y)>(i\circ c)(y)=x_1^0. $$
From Theorem \ref{thm.homotopy.alexandroff.spaces}, we conclude that $\textnormal{id}_\Delta$ is homotopic to $i\circ c$ since $\text{id}_\Delta\leq S\geq i\circ c$. Therefore, $\Delta$ is homotopy equivalent to $\{ x_1^0\}$.

\end{proo}

This result shows that the non-Hausdorff join can be interpreted as a way to ``kill'' homotopy when it is iterated infinitely many times because, in general, the finitely iterated non-Hausdorff join is not contractible unless we make the join with a contractible space (see, for example, Proposition \ref{prop_no_contractible}).  The following result is a trivial consequence of the definition of the non-Hausdorff join.

Let us recall that $Aut(X)$ denotes the group of self-homeomorphisms of $X$, it is immediate that $Aut(X\circledast Y)=Aut(X)\times Aut(Y)$. Moreover, $\mathcal{E}(X)$ denotes the group of self-homotopy equivalences of $X$. 

\begin{prop}\label{prop:automorphism_group} Let  $\{X_n \}_{n\in \mathbb{N}}$ be a sequence of non-empty finite spaces and consider $\Delta^{(n)}=X_{0}\circledast X_{1}\circledast \cdots \circledast X_{n} $ and $\Delta=X_{0}\circledast X_{1}\circledast \cdots \circledast X_{i} \circledast X_{i+1}\circledast\cdots$. Then $Aut(\Delta^{(n)})=Aut(X_0)\times Aut(X_1)\times \cdots \times Aut(X_n)$ and $Aut(\Delta)=\prod_{i\in \mathbb{N}} Aut(X_i)$.
\end{prop}

\begin{figure}[h!!]
\centering
\includegraphics[scale=1]{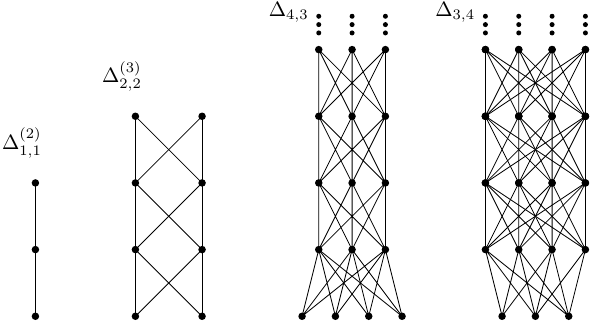}
\caption{From left to right, Hasse diagrams of $\Delta_{1,1}^{(2)}$, $\Delta_{2,2}^{(3)}$, $\Delta_{4,3}$ and $\Delta_{3,4}$.}\label{fig:Hasse diagram_delta_m_q}
\end{figure}

We now restrict our attention to the study of the non-Hausdorff join of the iterated join of antichains. Let $[q]$ denote the antichain of $q$ points. Suppose $m$, $q$ and $n$ are positive integer numbers, let $\Delta_{m,q}^{(n)}$ denote $[m]\circledast [q]\circledast \cdots^{n)} \circledast [q]$ and let $\Delta_{m,q}=[m]\circledast [q]\circledast [q]\circledast \cdots$. In Figure \ref{fig:Hasse diagram_delta_m_q}, we present the Hasse diagrams of $\Delta_{m,q}^{(n)}$ and $\Delta_{m,q}$ for some values of $m$, $q$ and $n$. 

\textbf{Notation:} Note that $\Delta_{m,q}^{(n)}$ denotes, with a slight abuse of notation, both the simplicial complex defined at the beginning of Section \ref{sect.qc} and the finite partially ordered set defined above. We prefer this notation to avoid an overload of symbols. We believe that the context will make it clear whether $\Delta_{m,q}^{(n)}$ refers to a simplicial complex or a partially ordered set.

\begin{prop}\label{prop_no_contractible}
    Let $m,q$ be positive integers and $n\in \mathbb{N}$. Then $\Delta_{m,q}^{(n)}$ is contractible if and only if $m$ or $q$ equals $1$. 
\end{prop}
\begin{proo}
It is clear that if $m,q>1$, then $\Delta_{m,q}^{(n)}$ does not have up or down beat points. For every $x\in \Delta_{m,q}^{(n)}$, we obtain that $U_{x}\setminus \{x \}=\emptyset$ (if $\textnormal{ht}(x)=0$) or $U_{x}\setminus \{x \}$ has $m$ or $q$ maximal points. The similar situation also holds when considering $L_x\setminus \{x \}$.

If $m=1$, then $\Delta_{m,q}^{(n)}$ has a minimum, which gives that it is contractible. If $q=1$, then $\Delta_{m,q}^{(n)}$ has a maximum and, consequently, contractible.
\end{proo}

\begin{rem} Note that we can obtain important examples of simplicial complexes from this construction by considering the order complex of $\Delta_{m,q}^{(n)}$. Let us enumerate some of them:
\begin{enumerate}
\item The simplicial complex $\mathcal{K}(\Delta_{1,1}^{(n)})$ is the simplex of dimension $n$.
\item The simplicial complex $\mathcal{K}(\Delta_{2,2}^{(n)})$ is the $n$-dimensional cross-polytope.
\item The simplicial complex $\mathcal{K}(\Delta_{m,q}^{(1)})$ has the homotopy type of the wedge sum of $(m-1)(q-1)$ circles. Consequently, every simplicial complex of dimension $1$ is homotopy equivalent to $\mathcal{K}(\Delta_{m,q}^{(1)})$ for some $m$ and $q$.
\end{enumerate}
\end{rem}

By Section \ref{sect.properties}, we know that $\mathcal{K}(\Delta_{m,q}^{(n)})$ is a partitionable simplicial complex. We can use Remark \ref{rem:interpretacion.exact.facing.particiones} to provide a combinatorial interpretation of an exact facing in terms of the finite barycentric subdivision. Specifically, an exact facing $\varphi$ for $\mathcal{K}(\Delta_{m,q}^{(n)})$ provides us with a partitioning into intervals of the finite barycentric subdivision of $\Delta_{m,q}^{(n)}$, where we are considering the empty set as a point. Let us clarify this with an illustrative example.

\begin{figure}[h!]
\centering
\includegraphics[scale=1]{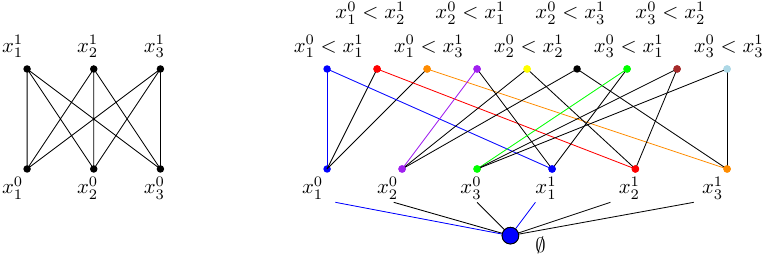}\caption{Hasse diagram of $\Delta_{3,3}^{(1)}$ on the left and Hasse diagram of $\mathcal{X}(\mathcal{K}(\Delta_{3,3}^{(1)}))$ with the empty set as a point, along with the partitioning induced by an exact facing, on the right.}\label{fig_ejemplo_exact_facing_particion_3_3}
\end{figure}

\begin{exam} Let us consider the finite space $\Delta_{3,3}^{(1)}=\{x_1^0,x_2^0,x_3^0,x_1^1,x_2^1,x_3^1 \}$ with $x_i^0<x_j^1$ for any $i,j=1,2,3$, whose Hasse diagram is depicted in Figure \ref{fig_ejemplo_exact_facing_particion_3_3}. Consider  the exact facing $\varphi:\textnormal{facets}(\mathcal{K}(\Delta_{3,3}^{(1)}))\rightarrow \mathcal{K}(\Delta_{3,3}^{(1)})$ defined by $\varphi(x_1^0<x_1^1)=\emptyset$, $\varphi(x_1^0<x_2^1)=x_2^1$, $\varphi(x_1^0<x_3^1)=x_3^1$, $\varphi(x_2^0<x_1^1)=x_2^0$, $\varphi(x_2^0<x_2^1)=x_2^0<x_2^1$, $\varphi(x_2^0<x_3^1)=x_2^0<x_3^1$, $\varphi(x_3^0<x_1^1)=x_3^0$, $\varphi(x_3^0<x_2^1)=x_3^0<x_2^1$, and $\varphi(x_3^0<x_3^1)=x_3^0<x_3^1$. Thus, the h-vector of $\mathcal{K}(\Delta_{3,3}^{(1)})$ is $(1,4,4)$, and $[\emptyset,x_1^0<x_1^1]$, $[x_2^1,x_1^0<x_2^1]$, $[x_3^1,x_1^0<x_3^1]$, $[x_2^0,x_2^0<x_1^1],[x_2^0<x_2^1,x_2^0<x_2^1]$, $[x_2^0<x_3^1,x_2^0<x_3^1]$, $[x_3^0,x_3^0<x_1^1]$, $[x_3^0<x_2^1,x_3^0<x_2^1]$, $[x_3^0<x_3^1,x_3^0<x_3^1]$ is a partitioning of $\mathcal{X}(\mathcal{K}(\Delta_{3,3}^{(1)}))$, with the empty set as a point, into disjoint intervals. We have represented each of these intervals with a different color in Figure \ref{fig_ejemplo_exact_facing_particion_3_3}.
\end{exam}

As an immediate consequence of results obtained in previous sections, we obtain the following result by considering the order complex of $\Delta_{m,q}^{(n)}$.
\begin{prop}\label{cor:Delta_contractible} Let $q,m$ and $n$ be positive integer numbers and $s=(m-1)(q-1)^n$. Then $\Delta_{m,q}^{(n)}$ is weak homotopy equivalent to the wedge sum of $s$ copies of the $n$-dimensional sphere.
\end{prop}

We will derive new Riordan matrices in the subsequent results, defined in terms of the topological or combinatorial invariants of the spaces $\Delta_{m,q}$.
\begin{theo}\label{thm:automorphism_group} Let  $m,q>1$ and $n$ be positive integers. Then 
\begin{enumerate}
\item $Aut(\Delta_{m,q}^{(n)})=\mathcal{E}(\Delta_{m,q}^{(n)})=S_m \times (\prod_{i=1}^n S_q)$.
\item $Aut(\Delta_{m,1}^{(n)})=S_m$ and $\mathcal{E}(\Delta_{m,1}^{(n)})$ is the trivial group.
\item $Aut(\Delta_{1,q}^{(n)})=\prod_{i=1}^n S_q$ and $\mathcal{E}(\Delta_{1,q}^{(n)})$ is the trivial group.
\item $Aut(\Delta_{1,1}^{(n)})=\mathcal{E}(\Delta_{1,1}^{(n)})$ is the trivial group.
\item $Aut(\Delta_{m,q})=S_m \times (\prod_{i=1}^\infty S_q)$ and $\mathcal{E}(\Delta_{m,q})$ is the trivial group.
\end{enumerate}
\end{theo}
\begin{proo}
The result follows from Proposition \ref{prop:automorphism_group} and the following facts: (1) The space $\Delta_{m,q}^{(n)}$ is a minimal finite space for every $n\geq 1$. (2) The space $\Delta_{m,1}^{(n)}$ is contractible (because it has a maximum) for every $n\geq 1$. (3) The space $\Delta_{1,q}^{(n)}$ is contractible (because it has a minimum) for every $n\geq 1$. (4) The space $\Delta_{1,1}^{(n)}$ is contractible for every $n\geq 1$. (5) The space $\Delta_{m,q}$ is contractible from Corollary \ref{cor:Delta_contractible}.
 \end{proo}

\begin{coro} Let $m$ and $q$ be positive integers and let $\bm C_{m,q} =(c_{i,j})_{i,j\geq0}$ be the matrix defined by $c_{i,j}=\textnormal{Card}(Aut(\Delta_{m,q}^{(i)}))$ if $i=j$ and $c_{i,j}=0$ otherwise (where \textnormal{Card} denotes the cardinal). Then
\begin{enumerate}
    \item $\bm{C}_{m,q}$ is the Riordan matrix $T\left(\frac{m!}{q!}\Bigm|\frac{1}{q!}\right)$ if $m,q\geq 1$.
    \item $\bm{C}_{m,q}$ is the identity matrix if $m=q=1$.
\end{enumerate}
\end{coro}
\begin{proo}
    It is an immediate consequence of Theorem \ref{thm:automorphism_group}.
    
\end{proo}

An important notion in the theory of partially ordered sets is the \textbf{order dimension}. Let $X$ be a partially ordered set. A linear extension $E$ of $P$ is a total ordering $\leq_1$ on $X$ such that $\text{id}:X\rightarrow E$ is order preserving. The dimension $\text{dim}(X)$ of a poset $X$ is the least $d$ for which there exists a family of linear extensions $E_1=(X,\leq_1),...,E_d=(X,\leq_d)$ such that 
$x\leq y$ if and only if $x\leq_1 y,\ldots,x\leq_d y$.

\begin{theo}\label{thm.poset.dimension}
Let $q,m$ and $n$ be positive integer numbers. Then $\textnormal{dim}(\Delta_{m,q}^{(n)})=\textnormal{dim}(\Delta_{m,q})=1$ if $m=q=1$, and $\textnormal{dim}(\Delta_{m,q}^{(n)})=\textnormal{dim}(\Delta_{m,q})=2$ otherwise. 
\end{theo}
\begin{proo}
If $m=q=1$, then $\Delta_{1,1}^{(n)}$ and $\Delta_{1,1}$ are totally ordered sets, and this means that $\text{dim}(\Delta_{m,q}^{(n)})=1$. Suppose now that $m$ or $q$ are distinct from $1$. We introduce notation: $\Delta_{m,q}^{(n)}=\{x^0_{i} \}_{i=1,...,m} \cup \{x_i^j \}_{\substack{i=1,...,q \\ j=1,...,n}}$, where $\text{ht}(x_i^j)=j$ for $j=0,...,n$. Consider $E_1=(\Delta_{m,q}^{(n)},\leq_1)$ as follows:
\begin{align*}
x_1^0<_1x_2^0<_1\cdots<_1x_m^0<_1 x_1^1<_1x_2^1<_1\cdots< x_q^1<_1x_1^2<_1x_2^2<_1\cdots \\\cdots<x_q^2<_1\cdots <_1x_1^n<_1x_2^n<_1\cdots <_1x_q^n.    
\end{align*}

Consider $E_2=(\Delta_{m,q}^{(n)},\leq_2)$ as follows:

\begin{align*}
    x_m^0<_2x_{m-1}^0<_2\cdots<_2x_1^0<_2 x_q^1<_2x_{q-1}^1<_2\cdots <x_1^1<_2x_q^2<_2x_{q-1}^2<_2\cdots\\\cdots <x_1^2<_2\cdots <_2x_q^n<_2x_{q-1}^n<_2\cdots <_2x_1^n.
\end{align*}

It is straightforward to verify that $x\leq y$ in $\Delta_{m,q}^{(n)}$ if and only if $x\leq_1 y$ and $x\leq_2 y$. The same argument also holds for $\Delta_{m,q}$.

\end{proo}

Let $m$ and $q$ be positive integers. Consider $\Delta_{m,q}$ and let $d_i = \text{dim}(\Delta_{m,q}^{(i)})$. By Theorem \ref{thm.poset.dimension}, the sequence $(d_i)_{i \in \mathbb{N}}$ is given by $d_i = 1$ for every $i$ when $m = q = 1$, and $d_i = 2$ otherwise. By storing this sequence of numbers in an infinite diagonal matrix, we obtain a Riordan matrix.

We now highlight how some previous Riordan matrices can also be defined from Alexandroff spaces using the order complex. Let us recall that $\bm F_{m,q}$ is the Riordan matrix obtained in Subsection \ref{sub:f-vector_matrix_h-vector}, while $\bm B_{m,q}$ is the Riordan matrix from Subsection \ref{subsect.betti}

\begin{theo} Let $\bm F'_{m,q}=(f_{n,k})_{n,k\geq0}$ be the infinite matrix such that the $n$-row of $\bm F'$ is just the $f$-vector of $\Delta^{(n)}_{m,q}$. Then $\bm F'_{m,q}=\bm F_{m,q}$. Moreover, the infinite matrix $\bm B_{m,q}'=(\tilde{\be'}_{n,k})_{n,k\geq0}$, where $\tilde{\be'}_{n,k}=\textnormal{rank}\  \tilde{H}_k(\Delta^{(n)}_{q,m})$, is the same as $\bm B_{m,q}$.
\end{theo}
\begin{proo}
    The first result is an immediate consequence of Remark \ref{rem.f.vectores.son.iguales.finitos.mccord}. The second result is evident from the results recalled in Subsection \ref{subsection.preliminares.posets}.
\end{proo}

To conclude this section, we try to follow the same steps to obtain a Riordan matrix from two numerical invariants, as it was done for the order dimension. From Corollary \ref{cor:Delta_contractible}, it is evident that the Euler characteristic of $\Delta_{m,q}^{(n)}$ is equal to $1+(-1)^{n}(m-1)(q-1)^{n}$, let us denote this number by $\chi_n$ and consider the matrix $(x_{i,j})_{i,j\geq 0}$ defined by $x_{i,j}=\chi_i$ if $i=j>0$  and $x_{i,j}=0$ otherwise. This matrix is not a Riordan matrix. On the contrary, if we consider the determinant of a finite space, then we obtain trivially that $\det(\Delta_{m,q}^{(n)})=(m-1)(q-1)^n$. Recall that the determinant of a finite space $X=\{x_1,...,x_n \}$ is provided by the absolute value of the determinant of the matrix $M_X=(x_{i,j})_{i,j=1,...,n}$ defined by $x_{i,j}=0$ if $x_i\leq x_j$ and $x_{i,j}=1$ otherwise. This number is a homotopy and simple homotopy invariant (see \cite{Chocano2} for more information). By storing these numbers in a matrix $\bm N_{m,q}=(d_{i,j})_{i,j\geq 0}$ as follows $d_{i,j}=\det(\Delta_{m,q}^{(n)})$ if $i=j=n$ and $d_{i,j}=0$ otherwise, we obtain a Riordan matrix. Particularly, we obtain the following result.

\begin{theo} Let $m$ and $q$ be positive integers. Then $\bm {N}_{m,q}=\bm{B}_{m,q}$.
\end{theo}

\subsection{Examples of iterated $q$-cones} \label{subsect.exam}


\subsubsection*{Case $m=q=1$. Cones and balls.}

In this case the Riordan matrix of the
$f$-vector of $\Delta_{1,1}^{(n)}$ is
\begin{small}
\begin{equation} \label{eq.F11}
\bm F_{11}=T\left(\frac{1}{1-x}\Big|1-x\right)=
\left(%
\begin{array}{cccccc}
  1 &  &  &  &  \\
  2 & 1 &  &  &  \\
  3 & 3 & 1 &  &  \\
  4 & 6 & 4 & 1 &  \\
  5 & 10 & 10 & 5 & 1 \\
  \vdots & \vdots & \vdots & \vdots & \vdots & \ddots \\
\end{array}%
\right)
\end{equation}
\end{small}

\noindent which is the matrix obtained from the Pascal triangle by deleting
the first row and the first column. So
\begin{equation}\label{eq.11fnk}
f_{n,k}=\binom{n+1}{k+1}, \qquad n,k\geq0.
\end{equation}
The corresponding polyhedra are, topologically, closed balls in
Euclidean spaces. Moreover
\[
\chi_{1,1}(x)=\frac{1}{1-x}=\sum_{n\geq0}x^n, \qquad \text{and}
\qquad \chi(\Delta_{1,1}^{(n)})=1 \quad \forall \ n\geq0
\]
which reflects, of course, the fact that the Euler characteristic
is always $1$ because all polyhedra in the sequence have the
homotopy type of a point.

\begin{figure}[h!]
\centering
\includegraphics[width=0.7\textwidth]{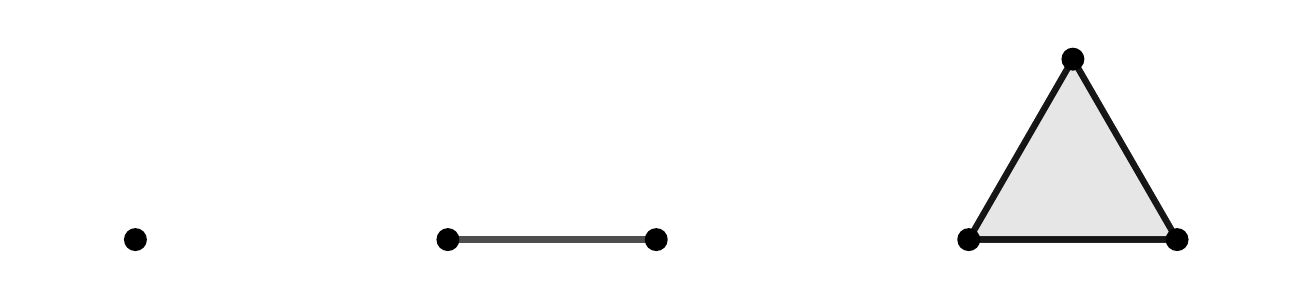}
\caption{The simplicial complexes $\Delta_{11}^{(0)}$, $\Delta_{11}^{(1)}$ and $\Delta_{11}^{(2)}$.}

\end{figure}

\subsubsection*{Case $m=q=2$. Suspensions and spheres.}

If $m=q=2$, the Riordan matrix we get is
\begin{small}
\begin{equation} \label{eq.F22}
T\left(\frac{1}{1-x}\Big|\frac{1-x}{2}\right)=
\left(%
\begin{array}{cccccc}
  2 &  &  &  &  \\
  4 & 4 &  &  &  \\
  6 & 12 & 8 &  &  \\
  8 & 24 & 32 & 16 &  \\
  10 & 40 & 80 & 80 & 32 \\
  \vdots & \vdots & \vdots & \vdots & \vdots & \ddots \\
\end{array}%
\right)
\end{equation}
\end{small}
or
\begin{equation} \label{eq.22fnk}
f_{n,k}=\binom{n+1}{k+1}2^{k+1}.
\end{equation}

\begin{figure}[h!]
\centering
\includegraphics[width=0.7\textwidth]{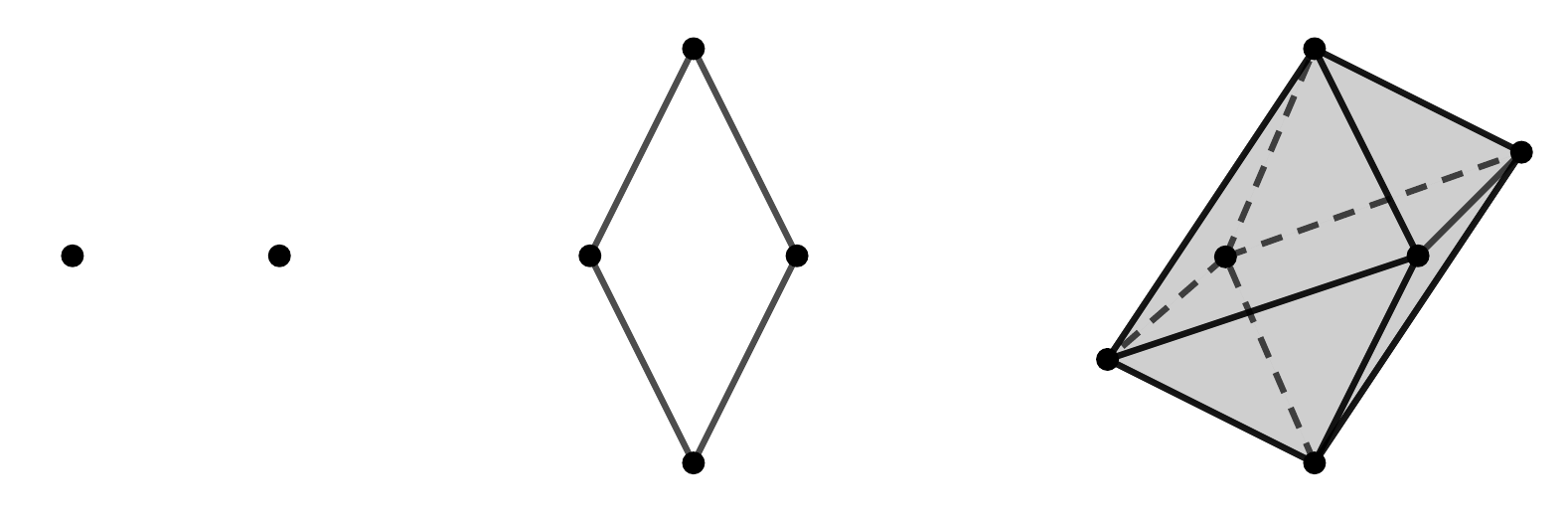}
\caption{The simplicial complexes $\Delta_{22}^{(0)}$, $\Delta_{22}^{(1)}$ and $\Delta_{22}^{(2)}$.}

\end{figure}

Now
$$\chi_{2,2}(x)=\frac{1}{1-x^2}=\sum_{n\geq0}(1+(-1)^n)x^n $$

\noindent So
$$f_{n0}-f_{n1}+f_{n2}-\ldots+(-1)^nf_{nn}=1+(-1)^n \quad \forall
\ n\geq0.$$

\noindent Note also that Corollary \ref{coro.meuler} gives us that
$$\frac{1}{2}f_{n0}-\frac{1}{4}f_{n1}+\frac{1}{8}f_{n2}-\ldots+(-1)^n\frac{1}{2^{n+1}}f_{nn}=1. $$

Let us remark that the two previous linear relations are linearly independent.

\subsubsection*{General case $m=1, q\geq 1$.}

 Note that, in this case, the simplicial complexes $\Delta^{(n)}_{1,q}$ are, all of them, contractible, and then, their Euler characteristic should be 1, which is just what we have obtained making the algebraic computation
$$\bm F_{1,q}(M_1)=\frac{1}{1-x}. $$

On the other hand, let us consider the general case $m=2,q\geq 1$. The simplicial complexes $\Delta_{2,q}^{(n)}$ are more complicated. As an application of Equation \eqref{eq.eulercarqc} we have
$$\chi_{2,q}(x)=\sum_{n\geq 0}(1+(1-q)^n)x^n $$

\noindent and 
$$\chi(\Delta_{2,q}^{(n)})=f_{n0}-f_{n1}+\ldots+(-1)^n f_{nn}=1+(1-q)^n.$$

\noindent However, in this case, the formula
$$\bm F_{2,q}(M_2)=\frac{1}{1-x} $$

\noindent yields that, for every $n$, 
$$\frac{1}{2}f_{n0}-\frac{1}{4}f_{n1}+\frac{1}{8}f_{n2}-\ldots+(-1)^n\frac{1}{m^{n-1}}f_{nn}=1. $$


\end{document}